\documentclass[10pt,a4paper]{article}

\usepackage[leqno]{amsmath}
\usepackage{graphicx,latexsym,euscript,makeidx,color}
\usepackage{amsmath,amsfonts,amssymb,amsthm,mathrsfs}
\usepackage{enumerate}

\usepackage[colorlinks,linkcolor=blue,anchorcolor=green,citecolor=red]{hyperref}

\usepackage{geometry}
\def\5n{\negthinspace \negthinspace \negthinspace \negthinspace \negthinspace }
\def\4n{\negthinspace \negthinspace \negthinspace \negthinspace }
\def\3n{\negthinspace \negthinspace \negthinspace }
\def\2n{\negthinspace \negthinspace }
\def\1n{\negthinspace }

\def\dbE{\mathbb{E}}     
\def\dbF{\mathbb{F}}   \def\cF{{\cal F}}  
     
\def\dbH{\mathbb{H}}

   \def\cL{{\cal L}}  
   \def\cM{{\cal M}}  
   \def\cN{{\cal N}}  
     
\def\dbP{\mathbb{P}}     
     
\def\dbR{\mathbb{R}} \def\sR{\mathscr{R}}    
\def\dbS{\mathbb{S}} \def\sS{\mathscr{S}}    
     
   \def\cU{{\cal U}}  
   \def\cV{{\cal V}}  
     
\def\dbX{\mathbb{X}}   \def\cX{{\cal X}}

\def\ss{\smallskip}                \def\lt{\left}
\def\ms{\medskip}                \def\rt{\right}
         \def\Ra{\mathop{\Rightarrow}}      \def\lan{\langle}
\def\ds{\displaystyle}           \def\ran{\rangle}
\def\no{\noindent}        \def\q{\quad}                      \def\llan{\left\langle}
\def\ns{\noalign{\ss}}    \def\qq{\qquad}                    \def\rran{\right\rangle}
    \def\hb{\hbox}                     \def\blan{\big\langle}
                   \def\bran{\big\rangle}
         \def\rf{\eqref}                    \def\Blan{\Big\langle}
  \def\deq{\triangleq}               \def\Bran{\Big\rangle}
 \def\ae{\hbox{\rm a.e.}}           \def\({\Big (}
\def\les{\leqslant}       \def\as{\hbox{\rm a.s.}}           \def\){\Big )}
\def\ges{\geqslant}          \def\[{\Big[}
\def\cl{\overline}           \def\]{\Big]}
\def\h{\widehat}          \def\tr{\hbox{\rm tr$\,$}}         \def\cd{\cdot}
\def\wt{\widetilde}              
           \def\cl{\overline}

\def\a{\alpha}        \def\G{\Gamma}   \def\g{\gamma}
\def\b{\beta}         \def\D{\Delta}   \def\d{\delta}
\def\z{\zeta}         \def\Th{\Theta}  
\def\e{\varepsilon}   \def\L{\Lambda}  \def\l{\lambda}
    \def\Si{\Sigma}  \def\si{\sigma}
           \def\F{\Phi}     
\def\n{\nu}           \def\Om{\Omega}  
          \def\f{\varphi}  \def\i{\infty}

\def\bde{\begin{definition}\label}    \def\ede{\end{definition}}
\def\bt{\begin{theorem}\label}        \def\et{\end{theorem}}
\def\bc{\begin{corollary}\label}      \def\ec{\end{corollary}}
\def\bl{\begin{lemma}\label}          \def\el{\end{lemma}}
\def\bp{\begin{proposition}\label}    \def\ep{\end{proposition}}
\def\bas{\begin{assumption}\label}    \def\eas{\end{assumption}}
\def\br{\begin{remark}\label}         \def\er{\end{remark}}
\def\bex{\begin{example}\label}       \def\ex{\end{example}}
\def\ba{\begin{array}}                \def\ea{\end{array}}
\def\be{\begin{equation}}
\def\bel{\begin{equation}\label}      \def\ee{\end{equation}}
\def\bea{\begin{eqnarray*}}           \def\eea{\end{eqnarray*}}

\newtheorem{theorem}{\indent Theorem}[section]
\newtheorem{definition}[theorem]{\indent Definition}
\newtheorem{proposition}[theorem]{\indent Proposition}
\newtheorem{corollary}[theorem]{\indent Corollary}
\newtheorem{lemma}[theorem]{\indent Lemma}
\newtheorem{remark}[theorem]{\indent Remark}
\newtheorem{example}[theorem]{\indent Example}

\newtheorem{assumption}[theorem]{\indent Assumption}

\makeatletter
   
   \@addtoreset{equation}{section}
\makeatother

\allowdisplaybreaks[4]

\begin{document}

\title{\bf Stochastic Linear Quadratic Optimal Control Problems in Infinite Horizon
}
\author{Jingrui Sun\thanks{Department of Mathematics, National University of Singapore, 119076,
Singapore (sjr@mail.ustc.edu.cn).}\,\, \ \ and \ \ Jiongmin Yong\thanks{Department of
Mathematics, University of Central Florida, Orlando, FL 32816, USA (jiongmin.yong@ucf.edu).
This author was partially supported by NSF grant DMS-1406776.}
}
\maketitle

\no {\bf Abstract:}
This paper is concerned with stochastic linear quadratic (LQ, for short) optimal control problems
in an infinite horizon with constant coefficients. It is proved that the non-emptiness of the
admissible control set for all initial state is equivalent to the $L^2$-stabilizability of the
control system, which in turn is equivalent to the existence of a positive solution to an algebraic
Riccati equation (ARE, for short). Different from the finite horizon case, it is shown that both
the open-loop and closed-loop solvabilities of the LQ problem are equivalent to the existence of a
{\it static stabilizing solution} to the associated generalized ARE. Moreover, any open-loop optimal
control admits a closed-loop representation. Finally, the one-dimensional case is worked out completely
to illustrate the developed theory.

\ms

\no {\bf Key words:}
stochastic linear quadratic optimal control, stabilizability, open-loop solvability, closed-loop solvability,
algebraic Riccati equation, static stabilizing solution, closed-loop representation

\ms

\no\bf AMS subject classifications. \rm 49N10, 49N35, 93D15, 93E20

\section{Introduction}

Let $(\Om,\cF,\dbF,\dbP)$ be a complete filtered probability space on which a standard one-dimensional
Brownian motion $W=\{W(t); 0\les t<\i\}$ is defined, where $\dbF=\{\cF_t\}_{t\ges0}$ is the natural
filtration of $W$ augmented by all the $\dbP$-null sets in $\cF$.
For a Euclidean space $\dbH$, let $L_\dbF^2(\dbH)$ denote the space of $\dbF$-progressively measurable
processes $\f:[0,\i)\times\Om\to\dbH$ with $\dbE\int_0^\i|\f(t)|^2dt<\i$.

\ms

Consider the following controlled linear stochastic differential equation (SDE, for short)
on the infinite time horizon $[0,\i)$:
\bel{state}\left\{\2n\ba{ll}
\ds dX(t)=\big[AX(t)+Bu(t)+b(t)\big]dt+\big[CX(t)+Du(t)+\si(t)\big]dW(t),\qq t\ges0, \\
\ns\ds X(0)= x,\ea\right.\ee
with quadratic cost functional
\bel{cost}J(x;u(\cd))\deq
\dbE\int_0^\i\left[\llan\begin{pmatrix}Q&S^\top\\S&R\end{pmatrix}
                        \begin{pmatrix}X(t)\\ u(t)\end{pmatrix},
                        \begin{pmatrix}X(t)\\ u(t)\end{pmatrix}\rran
+2\llan\begin{pmatrix}q(t)\\ \rho(t)\end{pmatrix},
       \begin{pmatrix}X(t)\\ u(t)   \end{pmatrix}\rran\right]dt.\ee
Here and throughout the paper, $A,C,Q\in\dbR^{n\times n}$, $B,D,S^\top\1n\in\dbR^{n\times m}$,
and $R\in\dbR^{m\times m}$ are given constant matrices, with $Q$ and $R$ being symmetric;
the superscript $\top$ denotes the transpose of matrices and vectors;
and $b(\cd), \si(\cd), q(\cd)\in L_\dbF^2(\dbR^n)$, $\rho(\cd)\in L_\dbF^2(\dbR^m)$.
In \rf{state}, $X(\cd)$, valued in $\dbR^n$, is called the {\it state process} with {\it initial state}
$x\in\dbR^n$, and $u(\cd)$, which belongs to $L_\dbF^2(\dbR^m)$, is called the {\it control process}.
Note that for $(x,u(\cd))\in\dbR^n\times L_\dbF^2(\dbR^m)$, the solution $X(\cd)\equiv X(\cd\,;x,u(\cd))$
of \rf{state} might merely be locally square-integrable, i.e.,
$$\dbE\int_0^T|X(t)|^2dt<\i,\q\forall\, 0<T<\i;\qq\dbE\int_0^\i|X(t)|^2dt=\i,$$
and the above cost functional $J(x; u(\cd))$ might not be defined. Therefore, we introduce the following:
$$\cU_{ad}(x)\deq\lt\{u(\cd)\in L_\dbF^2(\dbR^m)~\bigg|~\dbE\int_0^\i|X(t;x,u(\cd))|^2dt<\i\rt\},
\qq x\in\dbR^n.$$
Any element $u(\cd)\in\cU_{ad}(x)$ is called an {\it admissible control} associated with $x$.
Our linear quadratic (LQ, for short) optimal control problem can now be stated as follows.

\ms

\bf Problem (LQ). \rm For any given $x\in\dbR^n$, find a $u^*(\cd)\in\cU_{ad}(x)$ such that
\bel{opt}J(x;u^*(\cd))=\inf_{u(\cd)\in\cU_{ad}(x)}J(x;u(\cd))\deq V(x).\ee

Any $u^*(\cd)\in\cU_{ad}(x)$ satisfying \rf{opt} is called an {\it open-loop optimal control}
of Problem (LQ) for the initial state $x$; the corresponding $X^*(\cd)\equiv X(\cd\,;x,u^*(\cd))$
is called an {\it optimal state process}; and the function $V(\cd)$ is called the {\it value function}
of Problem (LQ). In the special case that $b(\cd), \si(\cd), q(\cd), \rho(\cd)=0$, we denote the
corresponding LQ problem, cost functional and value function by Problem (LQ)$^0$, $J^0(x;u(\cd))$
and $V^0(x)$, respectively.

\ms

The study of LQ problems (both in finite and infinite horizons, deterministic and stochastic)
has a long history that can be traced back to the pioneering works of Kalman \cite{Kalman 1960}
and Wonham \cite{Wonham 1968} (see also \cite{Davis 1977, Bensoussan 1982, Anderson-Moore 1989,
Yong-Zhou 1999} and the references therein). In the literature, take the LQ problem in a finite
horizon as an example, it is typically assumed that the cost functional has a (uniformly) positive
definite weighting matrix for the control term and positive semidefinite weighting matrices for
the state terms. In such a case, the cost functional is convex and coercive in the control variable,
and thus the LQ problem (in a finite horizon) has a unique open-loop optimal control which further
admits a state feedback representation via the solution to a differential Riccati equation.
The problem in infinite horizon is similar, with some additional stabilizability conditions.

\ms

In 1983, while studying LQ problems on Hilbert space, You \cite{You 1983} found that the
weighting matrices of state in the cost functional do not have to be positive semidefinite.
In 1998, Chen--Li--Zhou \cite{Chen-Li-Zhou 1998} further found that for stochastic LQ problems,
even the control weighting matrix does not have to be positve semidefinite. Since then,
extensive research efforts have been devoted to the indefinite LQ problems and to the solvability
of indefinite Riccati equations (see, for examples, \cite{Lim-Zhou 1999, Chen-Zhou 2000,
Chen-Yong 2001, Ait Rami-Moore-Zhou 2001, Hu-Zhou 2003, Qian-Zhou 2013}).

\ms

Recently, Sun--Li--Yong \cite{Sun-Li-Yong 2016} found that an indefinite LQ problem with
finite horizon might admit an open-loop optimal control (which could even be continuous)
for any initial pair, while the Riccati equation is not solvable. This leads to the
introduction of the open-loop and closed-loop solvabilities which are essentially different.
It was shown in \cite{Sun-Li-Yong 2016} that the closed-loop solvability is actually equivalent
to the existence of a {\it regular solution} to the Riccati equation (such a fact was firstly
revealed by Sun--Yong \cite{Sun-Yong 2014} in 2014 for two-person zero-sum differential games),
and that the {\it strongly regular solvability} of the Riccati equation is equivalent to the
{\it uniform convexity} of the cost functional. See Yong \cite{Yong 2013}, Sun \cite{Sun 2016},
and Li--Sun--Yong \cite{Li-Sun-Yong 2016} for relevant works of LQ problems involving mean-field.

\ms

Stochastic LQ problems in an infinite horizon and the associated algebraic Riccati equations
(AREs, for short) were treated in \cite{Rami-Zhou 2000, Rami-Zhou-Moore 2000} via the linear
matrix inequality and semidefinite programming techniques. This approach has been further
developed extensively in \cite{Yao-Zhang-Zhou 2001} by virtue of a duality analysis of
semidefinite programming. Along another line, Wu--Zhou \cite{Wu-Zhou 2001} introduced a new
frequency characteristic and applied the classical frequency domain approach to the LQ problems.
Recently, based on the work of Yong \cite{Yong 2013}, Huang--Li--Yong \cite{Huang-Li-Yong 2015}
studied a mean-field LQ optimal control problem on $[0,\i)$.

\ms

In the research on LQ problems with infinite horizon, however, there are some important
issues that have not been addressed: \\[-1.8em]
\begin{enumerate}[\indent\rm(a)]
\item The precise relation between the solvability of the problem and that of the ARE.
      In \cite{Rami-Zhou-Moore 2000}, it was shown that the solvability of an LQ problem
      is equivalent to the existence of a {\it stabilizing solution} to the ARE,
      but it is hard to verify whether a solution is stabilizing since certain selection
      of additional processes is involved.\\[-1.8em]
\item The closed-loop solvability. Because closed-loop controls are independent of
      the initial state and the future information, it is more meaningful and convenient
      to use closed-loop controls rather than open-loop controls. On the other hand,
      closed-loop solvability trivially implies open-loop solvability.\\[-1.8em]
\item The structure of admissible control sets. Unlike the finite horizon case, the
      structure of $\cU_{ad}(x)$ seems to be very complicated since it involves the state
      equation. In general, $\cU_{ad}(x)$ depends on $x$, and may even be empty for some
      $x$. Figuring out the structure of admissible control sets will give us a better
      understanding of the LQ problem.\\[-1.8em]
\end{enumerate}
This paper is to address the above issues. An interesting fact we find is that for
infinite-horizon LQ problems, open-loop and closed-loop solvabilities coincide. Such a fact,
as we mentioned earlier, does not hold in the finite horizon case. We shall show that the
solvability of the problem can be characterized by the existence of a {\it static stabilizing
solution} to the ARE, and that every open-loop optimal control admits a closed-loop representation.
It is shown that the $L^2$-stabilizability is not only sufficient, but also necessary, for the
non-emptiness of all admissible control sets. Moreover, we show that the $L^2$-stabilizability can
be verified by solving an ARE.

\ms

The rest of this paper is organized as follows.
In Section 2, we present some preliminary results that shall be needed later.
Section 3 aims to describe the structure of admissible control sets.
In Section 4, we introduce the notions of open-loop and closed-loop solvabilities
as well as the algebraic Riccati equation, and state the main result of the paper.
Section 5 is devoted to a special case of Problem (LQ), where system \rf{state} is $L^2$-stable.
In Section 6 we prove the main result for the general case.
Finally, to illustrate the results obtained, we completely solve the one-dimensional case in section 7.

\section{Preliminaries}

{\bf 2.1. Notation.}
Let us first introduce/recall the following notation that will be used below:
\begin{eqnarray*}
&&\hb{$\dbR^{n\times m}$ is the Euclidean space of all $n\times m$ real matrices;
$\dbR^n=\dbR^{n\times 1}$ and $\dbR=\dbR^1$}.\\
&&\hb{$\dbS^n$: the space of all symmetric $n\times n$ real matrices}.\\
&&\hb{$\dbS^n_+$: the subset of $\dbS^n$ which consists of positive definite matrices}.\\
&&\hb{$\cl{\dbS^n_+}$: the subset of $\dbS^n$ which consists of positive semidefinite matrices}.\\
&&\hb{$M^\top$: the transpose of a matrix $M$}.\\
&&\hb{$M^\dag$: the Moore--Penrose pseudoinverse of a matrix $M$}.\\
&&\hb{$\tr(M)$: the sum of diagonal elements of a square matrix $M$}.\\
&&\hb{$|M|\deq\sqrt{\tr(MM^\top)}$: the Frobenius nrom of a matrix $M$}.\\
&&\hb{$\sR(M)$: the range of a matrix or an operator $M$}.
\end{eqnarray*}
For $M,N\in\dbS^n$, we use the notation $M\ges N$ (respectively, $M>N$) to indicate
that $M-N$ is positive semidefinite (respectively, positive definite).
Recall that the inner product $\lan\cd\,,\cd\ran$ on a Euclidean space is given by
$\lan M,N\ran\mapsto\tr(M^\top N)$.
When there is no confusion, we shall use $\lan\cd\,,\cd\ran$ for inner products in
possibly different Hilbert spaces, and denote by $|\cd|$ the induced norm.
Let $T>0$ and $\dbH$ be a Euclidean space. We denote
\begin{eqnarray*}
&& L_\dbF^2(0,T;\dbH)=\bigg\{\f:[0,T]\times\Om\to\dbH~\big|~\f(\cd)\hb{ is
$\dbF$-progressively measurable, } \dbE\int^T_0|\f(t)|^2dt<\i\bigg\},\\
&&  L_\dbF^2(\dbH)=\lt\{\f:[0,\i)\times\Om\to\dbH~\big|~\f(\cd)\hb{ is
$\dbF$-progressively measurable, } \dbE\int^\i_0|\f(t)|^2dt<\i\rt\},\\
&&  \cX[0,T]=\lt\{X:[0,\i)\times\Om\to\dbR^n~\big|~X(\cd)\hb{ is $\dbF$-adapted, continuous, }
\dbE\lt(\sup_{0\les t\les T}|X(t)|^2\rt)<\i\rt\}, \\
&&  \cX_{loc}[0,\i)=\bigcap_{T>0}\cX[0,T],\qq
\cX[0,\i)=\lt\{X(\cd)\in\cX_{loc}[0,\i)~\bigg|~\dbE\int_0^\i|X(t)|^2dt<\i\rt\}.
\end{eqnarray*}
We define the inner product on $L_\dbF^2(\dbH)$ by
$$\lan\f,\psi\ran=\dbE\int^\i_0\lan\f(t),\psi(t)\ran dt$$
so that $L_\dbF^2(\dbH)$ becomes a Hilbert space.

\ms

{\bf 2.2. $L^2$-stability.}
For given matrices $A,C\in\dbR^{n\times n}$, we denote by $[A,C]$ the following uncontrolled linear system:
\bel{equ-AC}\left\{\2n\ba{ll}
\ds dX(t)=AX(t)dt+CX(t)dW(t),\qq t\ges0, \\
\ns\ds X(0)=x.\ea\right.\ee

\bde{def-1}\rm System $[A,C]$ is said to be $L^2$-{\it stable} if for any $x\in\dbR^n$,
the solution $X(\cd\,;x)$ of \rf{equ-AC} is in $\cX[0,\i)$.
\ede

The following result, which will be used frequently in this paper, provides a characterization of the
$L^2$-stability of $[A,C]$. For a proof, see \cite{Rami-Zhou 2000, Huang-Li-Yong 2015}.

\bl{AC-stable}\sl The system $[A,C]$ is $L^2$-stable if and only if
there exists a $P\in\dbS^n_+$ such that
$$PA+A^\top P+C^\top PC<0.$$
In this case, for any $\L\in\dbS^n$, the Lyapunov equation
$$PA+A^\top P+C^\top PC+\L=0$$
admits a unique solution $P\in\dbS^n$ given by
$$P=\dbE\int_0^\i\F(t)^\top\L\F(t)dt,$$
where $\F(\cd)$ is the solution to the following SDE for $\dbR^{n\times n}$-valued processes:
\bel{F}\left\{\2n\ba{ll}
\ds d\F(t)=A\F(t)dt+C\F(t)dW(t),\qq t\ges0, \\
\ns\ds \F(0)=I.\ea\right.\ee
\el

Next, we present a result concerning the square-integrability of the solution to the system
\bel{L2-int}\left\{\2n\ba{ll}
\ds dX(t)=\big[AX(t)+b(t)\big]dt+\big[CX(t)+\si(t)\big]dW(t),\qq t\ges0, \\
\ns\ds X(0)=x.\ea\right.\ee
For the proof the reader is referred to Proposition 2.4 in Sun--Yong--Zhang \cite{Sun-Yong-Zhang 2016}.

\bl{COCV1}\sl Suppose that $[A,C]$ is $L^2$-stable.
Then for any $b(\cd),\si(\cd)\in L^2_\dbF(\dbR^n)$ and any $x\in\dbR^n$,
the solution $X(\cd\,;x,b(\cd),\si(\cd))$ of \rf{L2-int} is in $\cX[0,\i)$.
Moreover, there exists a constant $K>0$, independent of $x$, $b(\cd)$ and $\si(\cd)$ such that
$$\dbE\int_0^\i|X(t;x,b(\cd),\si(\cd))|^2dt\les K\lt[|x|^2+\dbE\int_0^\i\(|b(t)|^2+|\si(t)|^2\)dt\rt].$$
\el

We now consider the following backward stochastic differential equation (BSDE, for short) on $[0,\i)$:
\bel{IBSDE}dY(t)=-\big[A^\top Y(t)+C^\top Z(t)+\f(t)\big]dt+Z(t)dW(t),\qq t\in[0,\i).\ee

\bde{}\rm An $L^2$-{\it stable adapted solution} of \rf{IBSDE} is a pair
$(Y(\cd),Z(\cd))\in\cX[0,\i)\times L_\dbF^2(\dbR^n)$ which satisfies the integral version of \rf{IBSDE}:
$$Y(t)=Y(0)-\int_0^t\big[A^\top Y(s)+C^\top Z(s)+\f(s)\big]ds+\int_0^tZ(s)dW(s),
\qq\forall t\in[0,\i),\q\as$$
\ede

The following result, found in \cite{Sun-Yong-Zhang 2016},
establishes the existence and uniqueness of solutions to \rf{IBSDE}.

\bl{COCV2}\sl Suppose that $[A,C]$ is $L^2$-stable. Then for any $\f(\cd)\in L^2_\dbF(\dbR^n)$,
BSDE \rf{IBSDE} admits a unique $L^2$-stable adapted solution $(Y(\cd),Z(\cd))$.
\el

{\bf 2.3. Pseudoinverse.} We recall some properties of the pseudoinverse \cite{Penrose 1955}.

\bl{Penrose}\sl {\rm(i)} For any $M\in\dbR^{m\times n}$, there exists a unique matrix
$M^\dag\in\dbR^{n\times m}$ such that
$$MM^\dag M=M,\q M^\dag MM^\dag=M^\dag,\q (MM^\dag)^\top=MM^\dag,\q (M^\dag M)^\top=M^\dag M.$$
In addition, if $M\in\dbS^n$, then $M^\dag\in\dbS^n$, and
$$MM^\dag=M^\dag M;\qq M\ges0\iff M^\dag\ges0.$$

{\rm(ii)} Let $L\in\dbR^{n\times k}$ and $N\in\dbR^{n\times m}$.
The matrix equation $NX=L$ has a solution if and only if
\bel{Penrose-1}NN^\dag L=L,\ee
in which case the general solution is given by
$$X=N^\dag L+(I-N^\dag N)Y,$$
where $Y\in\dbR^{m\times k}$ is arbitrary.
\el

The matrix $M^\dag$ above is called the {\it Moore--Penrose pseudoinverse} of $M$.

\br{remark-1}\rm (i) Clearly, condition \rf{Penrose-1} is equivalent to $\sR(L)\subseteq\sR(N)$.

\ms

(ii) It can be easily seen from Lemma \ref{Penrose} that if $N\in\dbS^n$ and $NX=L$,
then $X^\top NX=L^\top N^\dag L$.
\er

{\bf 2.4. A useful lemma.}
We conclude this section with a lemma that will be used frequently in what follows.

\bl{bl-useful lemma}\sl
Let $\h Q\in\dbS^n$, $\h R\in\dbS^m$, and $\h S\in\dbR^{m\times n}$ be given.
Suppose that for each $T>0$ the differential Riccati equation
$$\left\{\2n\ba{ll}
\ds\dot P(t;T)+P(t;T)A+A^\top P(t;T)+C^\top P(t;T)C+\h Q\\
\ns\ds\q-\lt[P(t;T)B+C^\top P(t;T)D+\h S^\top\rt]\lt[\h R+D^\top P(t;T)D\rt]^{-1}\\
\ns\ds\qq~\cd\lt[B^\top P(t;T)+D^\top P(t;T)C+\h S\rt]=0,\qq t\in[0,T],\\
\ns\ds P(T;T)=G,
\ea\right.$$
admits a solution $P(\cd\,;T)\in C([0,T];\dbS^n)$ such that
$$\h R+D^\top P(t;T)D>0,\qq\forall t\in[0,T].$$
If $P(0;T)$ converges to $P$ as $T\to\i$ and $\h R+D^\top PD$ is invertible, then
$$PA+A^\top P+C^\top PC+\h Q
-\lt(PB+C^\top PD+\h S^\top\rt)\lt(\h R+D^\top PD\rt)^{-1}\lt(B^\top P+D^\top PC+\h S\rt)=0.$$
\el

\begin{proof}
For fixed but arbitrary $0<T_1<T_2<\i$, we define
$$\left\{\2n\ba{llll}
\ds P_1(t)\4n&=\4n&\ds P(T_1-t;T_1),\q& 0\les t\les T_1,\\
\ns P_2(t)\4n&=\4n&\ds P(T_2-t;T_2),\q& 0\les t\les T_2.
\ea\right.$$
On the interval $[0,T_1]$, both $P_1(\cd)$ and $P_2(\cd)$ satisfy the following equation for $\Si(\cd)$:
\bel{Ric-Sigma}\left\{\2n\ba{ll}
\ds \dot{\Si}(t)-\Si(t)A-A^\top\Si(t)-C^\top\Si(t)C-\h Q\\
\ns\hphantom{\dot{\Si}(t)}
+\[\Si(t)B+C^\top\Si(t)D+\h S^\top\]\[\h R+D^\top\Si(t)D\]^{-1}\[B^\top\Si(t)+D^\top\Si(t)C+\h S\]=0,\\
\ns \Si(0)=G.\ea\right.\ee
Set $\Pi(t)=P_1(t)-P_2(t)$, $0\les t\les T_1$. Then $\Pi(\cd)$ satisfies
$$\left\{\2n\ba{ll}
\ds \dot{\Pi}(t)-\Pi(t)A-A^\top\Pi(t)-C^\top\Pi(t)C\\
\ns\hphantom{\dot{\Pi}(t)}
+\big[\Pi(t)B+C^\top\Pi(t)D\big]\L_1(t)^{-1}\big[B^\top P_1(t)+D^\top P_1(t)C+\h S\,\big]\\
\ns\hphantom{\dot{\Pi}(t)}
-\big[P_2(t)B+C^\top P_2(t)D+\h S^\top\big]\L_2(t)^{-1}D^\top\Pi(t)D\L_1(t)^{-1}
\big[B^\top P_1(t)+D^\top P_1(t)C+\h S\,\big]\\
\ns\hphantom{\dot{\Pi}(t)}
+\big[P_2(t)B+C^\top P_2(t)D+\h S^\top\big]\L_2(t)^{-1}\big[B^\top\Pi(t)+D^\top\Pi(t)C\big]=0,
\qq t\in[0,T_1],\\
\ns \Pi(0)=0,\ea\right.$$
where $\L_1(t)=\h R+D^\top P_1(t)D>0$ and $\L_2(t)=\h R+D^\top P_2(t)D>0$.
Since $P_i(t)$ and $\L_i(t)^{-1}$, $i=1,2$ are uniformly bounded due to their continuity,
we may apply Gronwall's inequality to conclude $\Pi(t)\equiv0$. This shows
$$ P_1(\cd)=P_2(\cd),\q\hb{on }~[0,T_1]. $$
Therefore, we may define a function $\Si(\cd):[0,\i)\to\dbS^n$ by the following:
$$ \Si(t)=P(T-t;T),\q\hb{if }~0\les t\les T. $$
Noting that $\Si(\cd)$ satisfies \rf{Ric-Sigma} on the whole interval $[0,\i)$,
we have for any $T>0$,
\begin{eqnarray*}
&&\Si(T+1)-\Si(T)\\
&&=\int_T^{T+1}\Big\{\Si(t)A+A^\top\Si(t)+C^\top\Si(t)C+\h Q\\
&&\hphantom{=\int_T^{T+1}\Big\{}
-\[\Si(t)B+C^\top\Si(t)D+\h S^\top\]\[\h R+D^\top\Si(t)D\]^{-1}
\[B^\top\Si(t)+D^\top\Si(t)C+\h S\]\Big\}dt.
\end{eqnarray*}
The desired result follows by letting $T\to\i$ in the above.
\end{proof}

\section{Admissible Control Sets and Stabilizability}

In this section, we will look into the admissible control sets. Since the state equation is linear,
for any given $x\in\dbR^n$, the set $\cU_{ad}(x)$ of admissible controls is either empty or a convex
set in $L^2_\dbF(\dbR^m)$. To investigate Problem (LQ), we should find conditions for the system so
that the set $\cU_{ad}(x)$ is at least non-empty and hopefully it admits an accessible characterization.
To this end, we denote by $[A,C;B,D]$ the following controlled system:
$$\left\{\2n\ba{ll}
\ds dX(t)=\big[AX(t)+Bu(t)\big]dt+\big[CX(t)+Du(t)\big]dW(t),\qq t\ges0, \\
\ns\ds X(0)=x,\ea\right.$$
and introduce the following definition.

\bde{def-2}\rm System $[A,C;B,D]$ is said to be {\it $L^2$-stabilizable} if
there exists a $\Th\in\dbR^{m\times n}$ such that $[A+B\Th,C+D\Th]$ is $L^2$-stable.
In this case, $\Th$ is called a {\it stabilizer} of $[A,C;B,D]$.
The set of all stabilizers of $[A,C;B,D]$ is denoted by $\sS\equiv\sS[A,C;B,D]$.
\ede

One observes that the $L^2$-stabilizability of $[A,C;B,D]$ implies the non-emptiness
of the admissible control set $\cU_{ad}(x)$ for all $x$.
Indeed, if $\Th$ is a stabilizer of $[A,C;B,D]$, by Lemma \ref{COCV1}, the solution
$X(\cd)$ of the following SDE is in $\cX[0,\i)$:
$$\left\{\2n\ba{ll}
\ds dX(t)=\big[(A+B\Th)X(t)+b(t)\big]dt+\big[(C+D\Th)X(t)+\si(t)\big]dW(t),\qq t\ges0,\\
\ns\ds X(0)= x.\ea\right.$$
Hence, $u(\cd)\deq\Th X(\cd)\in\cU_{ad}(x)$.

\ms

The following result shows that the $L^2$-stabilizability is not only sufficient,
but also necessary, for the non-emptiness of all admissible control sets.

\bt{Uad-kehua}\sl The following statements are equivalent:

\ms

{\rm(i)} $\cU_{ad}(x)\ne\emptyset$ for all $x\in\dbR^n$;

\ms

{\rm(ii)} $\sS[A,C;B,D]\ne\emptyset$;

\ms

{\rm(iii)} The following ARE admits a positive solution $P\in\dbS^n_+$:
\bel{11.30}PA+A^\top P+C^\top PC+I
-\big(PB+C^\top PD\big)\big(I+D^\top PD\big)^{-1}\big(B^\top P+D^\top PC\big)=0.\ee

If the above are satisfied and $P$ is a positive solution of \rf{11.30}, then
\bel{12.1}\G\deq-\big(I+D^\top PD\big)^{-1}\big(B^\top P+D^\top PC\big)\in\sS[A,C;B,D].\ee
\et

\begin{proof}
We have proved the implication (ii) $\Ra$ (i) above.
For the implication (iii) $\Ra$ (ii), we observe that if $P$ is a positive solution
of \rf{11.30} and $\G$ is defined by \rf{12.1}, then
$$P(A+B\G)+(A+B\G)^\top P+(C+D\G)^\top P(C+D\G)=-I-\G^\top\G<0.$$
Hence, by Lemma \ref{AC-stable} and Definition \ref{def-2}, $\G$ is stabilizer of $[A,C;B,D]$.

\ms

We next show that (i) $\Ra$ (iii).
By subtracting solutions of \rf{state} corresponding to $x$ and $0$,
we may assume without loss of generality that $b(\cd)=\si(\cd)=0$.
Let $e_1,\cdots,e_n$ be the standard basis for $\dbR^n$.
Take $u_i(\cd)\in\cU_{ad}(e_i)$, $i=1,\cdots,n$, and set
$$U(\cd)=(u_1(\cd),\cdots,u_n(\cd)).$$
Then $U(\cd)x\in\cU_{ad}(x)$ for all $x\in\dbR^n$. Consider the following cost functional:
$$\bar J(x;u(\cd))=\dbE\int_0^\i\[|X(t)|^2+|u(t)|^2\]dt.$$
Let $\dbX(\cd)\in L_\dbF^2(\dbR^{n\times n})$ be the solution to the following SDE for
$\dbR^{n\times n}$-valued processes:
$$\left\{\2n\ba{ll}
\ds d\dbX(t)=\big[A\dbX(t)+BU(t)\big]dt+\big[C\dbX(t)+DU(t)\big]dW(t),\qq t\ges0,\\
\ns \dbX(0)=I.\ea\right.$$
We have
\begin{eqnarray*}
&&\inf_{u(\cd)\in\cU_{ad}(x)}\bar J(x;u(\cd))\les\dbE\int_0^\i\[|\dbX(t)x|^2+|U(t)x|^2\]dt\\
&&=\Blan\(\dbE\int_0^\i\big[\dbX(t)^\top\dbX(t)+U(t)^\top U(t)\big]dt\)x,x\Bran\equiv\lan\L x,x\ran,
\qq\forall x\in\dbR^n.
\end{eqnarray*}
Now for any fixed $T>0$, let us consider the state equation
$$\left\{\2n\ba{ll}
\ds dX_T(t)=\big[AX_T(t)+Bu(t)\big]dt+\big[CX_T(t)+Du(t)\big]dW(t),\qq t\in[0,T],\\
\ns X_T(0)=x,\ea\right.$$
and the cost functional
$$ \bar J_T(x;u(\cd))=\dbE\int_0^T\[|X_T(t)|^2+|u(t)|^2\]dt.$$
It is standard that the following differential Riccati equation
$$\left\{\2n\ba{ll}
\ds\dot P(t;T)+P(t;T)A+A^\top P(t;T)+C^\top P(t;T)C+I\\
\ns\ds~-\big[P(t;T)B+C^\top P(t;T)D\big]\big[I+D^\top P(t;T)D\big]^{-1}
\big[B^\top P(t;T)+D^\top P(t;T)C\big]=0,\q t\in[0,T],\\
\ns\ds P(T;T)=0,
\ea\right.$$
admits a unique solution $P(\cd\,;T)\in C([0,T];\dbS^n_+)$ such that
$$\inf_{u(\cd)\in L_\dbF^2(0,T;\dbR^m)}\bar J_T(x;u(\cd))
=\lan P(0;T)x,x\ran,\qq\forall x\in\dbR^n.$$
It is clear that
$$\ba{llll}
\ds \bar J_T(x;u(\cd))\les \bar J(x;u(\cd)),
&\ds\q\forall x\in\dbR^n,&\ds \forall u(\cd)\in \cU_{ad}(x),&\ds \forall\, 0<T<\i,\\
\ns\ds \bar J_T(x;u(\cd))\les \bar J_{T^\prime}(x;u(\cd)),
&\ds\q\forall x\in\dbR^n,&\ds\forall u(\cd)\in L_\dbF^2(0,T^\prime;\dbR^m),&\ds\forall\, 0<T<T^\prime<\i.
\ea$$
Thus, one has
$$0<P(0;T)\les P(0;T^\prime)\les \L,\qq 0<T<T^\prime<\i.$$
This implies that $P(0;T)$ converges increasingly to some $P\in\dbS^n_+$ as $T\nearrow \i$.
By Lemma \ref{bl-useful lemma}, $P$ solves the ARE \rf{11.30}.
\end{proof}

Theorem \ref{Uad-kehua} provides a characterization of the non-emptiness of
all admissible control sets. The following result further gives an explicit
description of the admissible controls.

\bp{p-ad}\sl Suppose that $\Th\in\sS[A,C;B,D]$. Then for any $x\in\dbR^n$,
$$\cU_{ad}(x)=\big\{\Th X_\Th(\cd\,;x,v(\cd))+v(\cd)\,|\,v(\cd)\in L^2_\dbF(\dbR^m)\big\},$$
where $X_\Th(\cd\,;x,v(\cd))$ is the solution to the following SDE:
\bel{X_Th}\left\{\2n\ba{ll}
\ds dX_\Th(t)=\big[(A+B\Th)X_\Th(t)+Bv(t)+b(t)\big]dt\\
\ns\ds\hphantom{dX_\Th(t)=}+\big[(C+D\Th)X_\Th(t)+Dv(t)+\si(t)\big]dW(t),\qq t\ges0, \\
\ns\ds X_\Th(0)= x.\ea\right.\ee
\ep

\begin{proof}
Let $v(\cd)\in L^2_\dbF(\dbR^m)$ and $X_\Th(\cd)\equiv X_\Th(\cd\,;x,v(\cd))$ be the
solution of \rf{X_Th}. Since $[A+B\Th,C+D\Th]$ is $L^2$-stable, by Lemma \ref{COCV1},
$X_\Th(\cd)\in\cX[0,\i)$. Set
$$u(\cd)\deq\Th X_\Th(\cd)+v(\cd)\in L^2_\dbF(\dbR^m),$$
and let $X(\cd)\in\cX_{loc}[0,\i)$ be the solution to
\bel{p-ad-1}\left\{\2n\ba{ll}
\ds dX(t)=\big[AX(t)+Bu(t)+b(t)\big]dt+\big[CX(t)+Du(t)+\si(t)\big]dW(t),\qq t\ges 0,\\
\ns X(0)=x.\ea\right.\ee
By uniqueness of solutions, we have $X(\cd)=X_\Th(\cd)\in\cX[0,\i)$,
and therefore $u(\cd)\in\cU_{ad}(x)$.

\ms

On the other hand, suppose $u(\cd)\in\cU_{ad}(x)$.
Let $X(\cd)\in\cX[0,\i)$ be the solution of \rf{p-ad-1} and set
$$v(\cd)\deq u(\cd)-\Th X(\cd)\in L^2_\dbF(\dbR^m).$$
Again by uniqueness of solutions, we see that $X(\cd)$ coincides with the solution
$X_\Th(\cd)$ of \rf{X_Th}.
Thus, $u(\cd)$ admits a representation of the form $\Th X_\Th(\cd\,;x,v(\cd))+v(\cd)$.
\end{proof}

To conclude this section, let us look at the case $n=1$, i.e., the state variable is one-dimensional.
First, we present the following lemma.

\bl{S-A}\sl Let $n=1$. If the system $[A,C;B,D]$ is not $L^2$-stabilizable, then
\bel{S-A-1}\begin{pmatrix}2A+C^2&B+CD\\B^\top\1n+CD^\top&D^\top D\end{pmatrix}\ges0.\ee
\el

\begin{proof}
If $[A,C;B,D]$ is not $L^2$-stabilizable, then by Definition \ref{def-2}
and Lemma \ref{AC-stable}, we have
$$2(A+B\Th)+(C+D\Th)^2\ges0,\qq\forall \Th\in\dbR^m.$$
Since for any nonzero $x\in\dbR$ and any $y\in\dbR^m$ one can find a $\Th\in\dbR^m$
such that $y=\Th x$, we have
\begin{eqnarray*}
&& \begin{pmatrix}x&\1ny^\top\end{pmatrix}
   \begin{pmatrix}2A+C^2&B+CD\\B^\top\1n+CD^\top&D^\top D\end{pmatrix}
   \begin{pmatrix}x\\y\end{pmatrix}\\
&& =x\begin{pmatrix}1&\1n\Th^\top\end{pmatrix}
     \begin{pmatrix}2A+C^2&B+CD\\B^\top\1n+CD^\top&D^\top D\end{pmatrix}
     \begin{pmatrix}1\\ \Th\end{pmatrix}x\\
&& =\big[2(A+B\Th)+(C+D\Th)^2\big]x^2\ges0,\qq\forall x\neq0,~y\in\dbR^m,
\end{eqnarray*}
and the result follows.
\end{proof}

For the moment let us assume that $b(\cd)=\si(\cd)=0$ and introduce the following space:
$$\cV=\big\{u(\cd)\in L_\dbF^2(\dbR^m)\,|\,Bu(\cd)=Du(\cd)=0,~\ae~\as\big\}.$$
Obviously, $0\in\cV\subseteq\cU_{ad}(0)$, and hence $\cU_{ad}(0)$ is non-empty. In fact,
when $\sS[A,C;B,D]=\emptyset$, $\cU_{ad}(0)$ coincides with $\cV$, and is the only non-empty
admissible control set. More precisely, we have the following result.

\bt{ad-1}\sl Let $n=1$, and suppose that $b(\cd)=\si(\cd)=0$.
Then exactly one of the following holds: \\[-1.6em]
\begin{enumerate}[\indent\rm(i)]
\item $\cU_{ad}(0)=\cV$ and $\cU_{ad}(x)=\emptyset$ for all $x\ne0$. \\[-1.6em]
\item The system $[A,C;B,D]$ is $L^2$-stabilizable.
\end{enumerate}\et

\begin{proof}
We prove it by contradiction. From Theorem \ref{Uad-kehua} we see that (i) and (ii)
cannot hold simultaneously. Now suppose that neither (i) nor (ii) holds.
Then either $\cU_{ad}(0)\backslash\cV\ne\emptyset$ or else $\cU_{ad}(x)\ne\emptyset$ for some $x\ne0$,
and \rf{S-A-1} holds by Lemma \ref{S-A}.
If there exists a $u(\cd)\in\cU_{ad}(0)\backslash\cV$, then with $X_0(\cd)$ denoting the solution of
\rf{state} corresponding to the initial state $x=0$ and the admissible control $u(\cd)$, we have
\begin{eqnarray*}
\dbE|X_0(t)|^2 \4n&=\4n& \dbE\int_0^t\Big\{2[AX_0(s)+Bu(s)]X_0(s)+|CX_0(s)+Du(s)|^2\Big\}ds\\
\4n&=\4n& \dbE\int_0^t\llan
\begin{pmatrix}2A+C^2 & B+CD\\B^\top\1n+CD^\top & D^\top D\end{pmatrix}
\begin{pmatrix}X_0(s)\\ u(s)\end{pmatrix},
\begin{pmatrix}X_0(s)\\ u(s)\end{pmatrix}\rran ds\ges0.
\end{eqnarray*}
Since \rf{S-A-1} holds and $X_0(\cd)\in\cX[0,\i)$ (and hence $\lim_{t\to\i}\dbE|X_0(t)|^2=0$),
the integrand in the above must vanish for all $s\ges0$. It turns out that $X_0(\cd)\equiv0$, and hence
$$Bu(\cd)=0, \qq Du(\cd)=0,\qq\ae~\as$$
which is a contradiction. Now if $\cU_{ad}(x)\ne\emptyset$ for some $x\ne0$, take $v(\cd)\in\cU_{ad}(x)$
and let $X(\cd)$ be the solution of \rf{state} corresponding to $x$ and $v(\cd)$. Then, using \rf{S-A-1},
we have for any $t\ges0$,
$$\dbE|X(t)|^2-|x|^2=\dbE\int_0^t\llan
\begin{pmatrix}2A+C^2 & B+CD\\B^\top\1n+CD^\top & D^\top D\end{pmatrix}
\begin{pmatrix}X(s)\\ v(s)\end{pmatrix},
\begin{pmatrix}X(s)\\ v(s)\end{pmatrix}\rran ds\ges0,$$
which is impossible since $\lim_{t\to\i}\dbE|X(t)|^2=0$. This completes the proof.
\end{proof}

For the case $b(\cd)\ne0$ or $\si(\cd)\ne0$, we have the following result,
which is a consequence of Theorem \ref{ad-1}.

\bt{ad-2}\sl Let $n=1$, and suppose that $b(\cd)\ne0$ or $\si(\cd)\ne0$.
Then exactly one of the following holds: \\[-1.6em]
\begin{enumerate}[\indent\rm(i)]
\item $\cU_{ad}(x)=\emptyset$ for all $x\in\dbR^n$. \\[-1.6em]
\item There is only one $x\in\dbR^n$ for which the admissible control set
$\cU_{ad}(x)\ne\emptyset$. In this case,
$$u(\cd)-v(\cd)\in\cV,\qq\forall u(\cd),v(\cd)\in\cU_{ad}(x).$$ \\[-3.6em]
\item The system $[A,C;B,D]$ is $L^2$-stabilizable.
\end{enumerate}\et

\begin{proof}
Clearly, any two of the statements (i)--(iii) cannot hold simultaneously.
To show the above, let us assume that neither (i) nor (ii) holds. Then
$$\cU_{ad}(x_1)\ne\emptyset,\qq\cU_{ad}(x_2)\ne\emptyset$$
for some $x_1\ne x_2$. Take $u_i(\cd)\in\cU_{ad}(x_i)$, $i=1,2$, and let $X_i(\cd)$ be the
solution of \rf{state} corresponding to the initial state $x_i$ and the admissible control
$u_i(\cd)$. Then with $x=x_1-x_2$ and $u(\cd)=u_1(\cd)-u_2(\cd)$, the process
$X(\cd)\deq X_1(\cd)-X_2(\cd)$ is in $\cX[0,\i)$ and solves
$$\left\{\2n\ba{ll}
\ds dX(t)=\big[AX(t)+Bu(t)\big]dt+\big[CX(t)+Du(t)\big]dW(t),\qq t\ges0, \\
\ns X(0)=x.\ea\right.$$
Thus, by Theorem \ref{ad-1}, the system $[A,C;B,D]$ is $L^2$-stabilizable.

\ms

Now suppose that there is only one $x\in\dbR^n$ such that $\cU_{ad}(x)\ne\emptyset$.
The same argument as before shows that for any $u(\cd),v(\cd)\in\cU_{ad}(x)$,
the solution $X_0(\cd)$ of
$$\left\{\2n\ba{ll}
\ds dX_0(t)=\big\{AX_0(t)+B[u(t)-v(t)]\big\}dt+\big\{CX_0(t)+D[u(t)-v(t)]\big\}dW(t),\qq t\ges0,\\
\ns\ds X_0(0)=0,\ea\right.$$
is in $\cX[0,\i)$. Since $\sS[A,C;B,D]=\emptyset$ in this situation,
we have $u(\cd)-v(\cd)\in\cV$ by Theorem \ref{ad-1}.
\end{proof}

\section{Solvabilities of Problem (LQ) and Generalized AREs}

Let us return to Problem (LQ). According to Theorem \ref{Uad-kehua},
when $[A,C;B,D]$ is not $L^2$-stabilizable, Problem (LQ) becomes ill-posed.
Because of this, we shall impose the following assumption in the rest of the paper:

\ms

{\bf(H1)} System $[A,C;B,D]$ is $L^2$-stabilizable, i.e., $\sS[A,C;B,D]\ne\emptyset$.

\ms

Now, we introduce the following definition.

\bde{def-3}\rm (i) An element $u^*(\cd)\in\cU_{ad}(x)$ is called an {\it open-loop optimal control}
of Problem (LQ) for the initial state $x\in\dbR^n$ if
\bel{open-opt}J(x;u^*(\cd))\les J(x;u(\cd)),\qq\forall u(\cd)\in\cU_{ad}(x).\ee
If an open-loop optimal control (uniquely) exists for $x\in\dbR^n$, Problem (LQ) is said to be
({\it uniquely}) {\it open-loop solvable at $x$}. Problem (LQ) is said to be ({\it uniquely})
{\it open-loop solvable} if it is (uniquely) open-loop solvable at all $x\in\dbR^n$.

\ms

(ii) A pair $(\Th^*,v^*(\cd))\in\dbR^{m\times n}\times L^2_\dbF(\dbR^m)$ is called a
{\it closed-loop optimal strategy} of Problem (LQ) if $\Th^*\in\sS[A,C;B,D]$ and
\bel{closed-opt} J(x;\Th^*X^*(\cd)+v^*(\cd))\les J(x;u(\cd)),
\qq\forall (x,u(\cd))\in\dbR^n\times\cU_{ad}(x),\ee
where $X^*(\cd)$ is the solution to the following closed-loop system:
\bel{close-X}\left\{\2n\ba{ll}
\ds dX^*(t)=\big[(A+B\Th^*)X^*(t)+Bv^*(t)+b(t)\big]dt\\
\ns\ds\hphantom{dX^*(t)=}+\big[(C+D\Th^*)X^*(t)+Dv^*(t)+\si(t)\big]dW(t), \qq t\ges0,\\
\ns\ds X^*(0)=x.\ea\right.\ee
If a closed-loop optimal strategy (uniquely) exists, Problem (LQ) is said to be ({\it uniquely})
{\it closed-loop solvable}.
\ede

It is worth pointing out that, in general, the admissible control sets $\cU_{ad}(x)$ are different
for different $x$, and an open-loop optimal control depends on the initial state $x\in\dbR^n$,
whereas a closed-loop optimal strategy is required to be independent of $x$.
From \rf{closed-opt}, one sees that the outcome $u^*(\cd)\equiv\Th^*X^*(\cd)+v^*(\cd)$ of
a closed-loop optimal strategy $(\Th^*,v^*(\cd))$ is an open-loop optimal control for the
initial state $X^*(0)$. Hence, closed-loop solvability implies open-loop solvability.
For LQ optimal control problems in finite horizon, the same is true, and open-loop solvability
does not necessarily imply closed-loop solvability (see \cite{Sun-Li-Yong 2016}).
However, for our Problem (LQ) (in an infinite horizon), as we shall prove later,
the open-loop and closed-loop solvabilities are equivalent, and both are equivalent to the existence
of a {\it static stabilizing solution} to a {\it generalized algebraic Riccati equation}
which we are going to introduce below.

\bde{def-ARE}\rm The following constrained nonlinear algebraic equation
\bel{ARE}\left\{\2n\ba{ll}
\ds PA+A^\top P+C^\top PC+Q\\
\ns\ds\hphantom{PA}
-\big(PB+C^\top PD+S^\top\big)\big(R+D^\top PD\big)^\dag\big(B^\top P+D^\top PC+S\big)=0,\\
\ns\ds\sR\big(B^\top P+D^\top PC+S\big)\subseteq\sR\big(R+D^\top PD\big),\\
\ns\ds R+D^\top PD\ges0 \ea\right.\ee
with the unknown $P\in\dbS^n$, is called a {\it generalized algebraic Riccati equation}.
A solution $P$ of \rf{ARE} is said to be {\it static stabilizing} if there exists a
$\Pi\in\dbR^{m\times n}$ such that
$$-\big(R+D^\top PD\big)^\dag\big(B^\top P+D^\top PC+S\big)
+\big[I-\big(R+D^\top PD\big)^\dag\big(R+D^\top PD\big)\big]\Pi\in\sS[A,C;B,D].$$
\ede

For notational simplicity, we shall write hereafter
$$\cM(P)=PA+A^\top P+C^\top PC+Q, \q~\cL(P)=PB+C^\top PD+S^\top, \q~\cN(P)=R+D^\top PD.$$

Now we state the main result of this paper.

\bt{main}\sl Let {\rm(H1)} hold. Then the following statements are equivalent:

\ms

{\rm(i)} Problem {\rm(LQ)} is open-loop solvable;

\ms

{\rm(ii)} Problem {\rm(LQ)} is closed-loop solvable;

\ms

{\rm(iii)} ARE \rf{ARE} admits a static stabilizing solution $P\in\dbS^n$, and the BSDE
\bel{eta-zeta}\ba{ll}
\ds d\eta=-\Big\{\big[A^\top\1n-\cL(P)\cN(P)^\dag B^\top\big]\eta+\big[C^\top\1n-\cL(P)\cN(P)^\dag D^\top\big]\z\\
\ns\ds\hphantom{d\eta=-\Big\{}
+\big[C^\top\1n-\cL(P)\cN(P)^\dag D^\top\big]P\si-\cL(P)\cN(P)^\dag\rho+Pb+q\Big\}dt+\z dW, \qq t\ges0,
\ea\ee
admits an $L^2$-stable adapted solution $(\eta(\cd),\z(\cd))$ such that
\bel{range-ez}B^\top\eta(t)+D^\top\z(t)+D^\top P\si(t)+\rho(t)\in\sR\big(\cN(P)\big),\q\ae~t\in[0,\i),~~\as\ee

In the above case, any closed-loop optimal strategy $(\Th^*,v^*(\cd))$ is given by
\bel{close-rep}\left\{\2n\ba{ll}
\ds \Th^*=-\cN(P)^\dag\cL(P)^\top+\big[I-\cN(P)^\dag\cN(P)\big]\Pi,\\
\ns v^*(\cd)=-\cN(P)^\dag\big[B^\top\eta(\cd)+D^\top\z(\cd)+D^\top P\si(\cd)+\rho(\cd)\big]
+\big[I-\cN(P)^\dag\cN(P)\big]\n(\cd),
\ea\right.\ee
where $\Pi\in\dbR^{m\times n}$ is chosen so that $\Th^*\in\sS[A,C;B,D]$, and $\n(\cd)\in L_\dbF^2(\dbR^m)$;
every open-loop optimal control $u^*(\cd)$ for the initial state $x$ admits a closed-loop representation:
\bel{u-star}u^*(\cd)=\Th^*X^*(\cd)+v^*(\cd),\ee
where $(\Th^*,v^*(\cd))$ is a closed-loop optimal strategy of Problem {\rm(LQ)} and $X^*(\cd)$ is the
corresponding solution to the closed-loop system \rf{close-X}. Further, the value function admits the
following representation:
\bel{Value}\ba{ll}
\ds V(x)=\lan Px,x\ran+\dbE\bigg\{2\lan\eta(0),x\ran
+\int_0^\i\[\lan P\si,\si\ran+2\lan\eta,b\ran+2\lan\z,\si\ran\\
\ns\ds\hphantom{V(x)=\,}
-\blan\cN(P)^\dag\big(B^\top\eta+D^\top\z+D^\top P\si+\rho\big),
B^\top\eta+D^\top\z+D^\top P\si+\rho\bran\]dt\bigg\}.\ea\ee
\et

The proof will be given in the subsequent sections. We make some observations here.
Suppose that ARE \rf{ARE} admits a static stabilizing solution $P\in\dbS^n$.
Then one can choose a matrix $\Pi\in\dbR^{m\times n}$ such that
$$\Th=-\cN(P)^\dag\cL(P)^\top+\big[I-\cN(P)^\dag\cN(P)\big]\Pi\in\sS[A,C;B,D].$$
If $(\eta(\cd),\z(\cd))$ is an $L^2$-stable adapted solution of \rf{eta-zeta} satisfying
\rf{range-ez}, then it follows easily from Lemma \ref{Penrose} that
$$\big[\Th^\top\1n+\cL(P)\cN(P)^\dag\big]\big[B^\top\eta(t)+D^\top\z(t)+D^\top P\si(t)+\rho(t)\big]=0.$$
Thus,
\begin{eqnarray*}
d\eta \4n&=\4n& -\,\Big\{A^\top\eta+C^\top\z+C^\top P\si+Pb+q
-\cL(P)\cN(P)^\dag\big[B^\top\eta+D^\top\z+D^\top P\si+\rho\big]\Big\}dt+\z dW\\
\4n&=\4n& -\,\Big\{A^\top\eta+C^\top\z+C^\top P\si+Pb+q
+\Th^\top\big[B^\top\eta+D^\top\z+D^\top P\si+\rho\big]\Big\}dt+\z dW\\
\4n&=\4n& -\,\Big\{(A+B\Th)^\top\eta+(C+D\Th)^\top\z+(C+D\Th)^\top P\si+Pb+q+\Th^\top\rho\Big\}dt+\z dW.
\end{eqnarray*}
Since $[A+B\Th,C+D\Th]$ is $L^2$-stable, we conclude from Lemma \ref{COCV2} that
the $L^2$-stable adapted solution of \rf{eta-zeta} satisfying \rf{range-ez} is unique.
In particular, when $b(\cd), \si(\cd), q(\cd), \rho(\cd)=0$, $(\eta(\cd),\z(\cd))\equiv(0,0)$
is such a solution. This leads to the following result.

\bc{main-bc}\sl Let {\rm(H1)} hold. Then the following statements are equivalent: \\[-1.6em]
\begin{enumerate}[\indent\rm(i)]
\item Problem {\rm(LQ)$^0$} is open-loop solvable; \\[-1.6em]
\item Problem {\rm(LQ)$^0$} is closed-loop solvable; \\[-1.6em]
\item ARE \rf{ARE} admits a static stabilizing solution $P\in\dbS^n$.\\[-1.6em]
\end{enumerate}

In the above case, all closed-loop optimal strategies $(\Th^*,v^*(\cd))$ are given by
\bel{LQ0-1}\Th^*=-\cN(P)^\dag\cL(P)^\top+\big[I-\cN(P)^\dag\cN(P)\big]\Pi,\qq
v^*(\cd)=\big[I-\cN(P)^\dag\cN(P)\big]\n(\cd),\ee
where $\Pi\in\dbR^{m\times n}$ is chosen so that $\Th^*\in\sS[A,C;B,D]$, and $\n(\cd)\in L_\dbF^2(\dbR^m)$;
the value function is given by
\bel{LQ0-2}V(x)=\lan Px,x\ran,\qq x\in\dbR^n.\ee
\ec

\begin{remark}\rm
From \rf{LQ0-2}, it is easily seen that ARE \rf{ARE} admits at most one static stabilizing solution.
\end{remark}

To conclude this section, we give an informal explanation why the open-loop solvability and
the closed-loop solvability coincide for LQ optimal control problems in infinite horizon.
We take Problem {\rm(LQ)$^0$} for example. In this case the state equation and the cost
functional respectively become
\begin{eqnarray}
&\label{16July18-15:30}\left\{\2n\ba{ll}
\ds dX(t)=\big[AX(t)+Bu(t)\big]dt+\big[CX(t)+Du(t)\big]dW(t),\qq t\ges0, \\
\ns X(0)= x,\ea\right.&\\
\ns&\label{}\ds
J(x;u(\cd))=\dbE\int_0^\i\llan\begin{pmatrix}Q&S^\top\\S&R\end{pmatrix}
                              \begin{pmatrix}X(t)\\ u(t)\end{pmatrix},
                              \begin{pmatrix}X(t)\\ u(t)\end{pmatrix}\rran dt.&
\end{eqnarray}
Suppose Problem (LQ)$^0$ is open-loop solvable. Then there exists an $\dbR^{m\times n}$-valued
process $U^*(\cd)$ such that for any initial distribution $\xi$, $U^*(\cd)\xi$ is the
open-loop optimal control corresponding to the linear quadratic structure of Problem (LQ)$^0$.
Let $X^*(\cd)$ be the solution of \rf{16July18-15:30} corresponding to the initial state $x$ and
the open-loop optimal control $u^*(\cd)\equiv U^*(\cd)x$. For any $s\ges0$, we may consider
Problem (LQ)$^0$ with initial time being $s$. That is to say, the state equation reads
$$\left\{\2n\ba{ll}
\ds dX(t)=\big[AX(t)+Bu(t)\big]dt+\big[CX(t)+Du(t)\big]dW(t),\qq t\ges s, \\
\ns\ds X(s)=x,\ea\right.$$
and the cost functional reads
$$J(s,x;u(\cd))=\dbE\int_s^\i\llan\begin{pmatrix}Q&S^\top\\S&R\end{pmatrix}
                                  \begin{pmatrix}X(t)\\ u(t)\end{pmatrix},
                                  \begin{pmatrix}X(t)\\ u(t)\end{pmatrix}\rran dt.$$
We denote this problem by Problem (LQ)$_s^0$. Since the matrices in the state equation and in the cost
functional are time-invariant and the time horizon is infinite, Problem (LQ)$^0$ and Problem (LQ)$_s^0$
can be regarded as the same problem. Thus, $U^*(\cd)X^*(s)$ is the open-loop optimal control of
Problem (LQ)$_s^0$ for the initial distribution $X^*(s)$.
On the other hand, by the dynamic programming principle, $U^*(s+\cd)x$ is also an open-loop optimal control
of Problem (LQ)$_s^0$ for the initial distribution $X^*(s)$. Therefore, one should have
$$U^*(s+t)x=U^*(t)X^*(s)\qq\forall s,t\ges0.$$
In particular, taking $t=0$, we have
$$u^*(s)=U^*(s)x=U^*(0)X^*(s),\qq s\ges0,$$
which should imply the closed-loop solvability of Problem (LQ)$^0$.
The key fact employed here is that, for any $s\ges0$, Problem (LQ)$_s^0$ is the same as Problem (LQ)$^0$.
We call such a property the {\it stationarity} of Problem (LQ)$^0$.

\section{The $L^2$-Stable Case}

In this section, we shall prove Theorem \ref{main} for the case $0\in\sS[A,C;B,D]$,
that is, $[A,C]$ is $L^2$-stable. Recall from Proposition \ref{p-ad} that in this case,
$$\cU_{ad}(x)=L_\dbF^2(\dbR^m),\qq\forall x\in\dbR^n.$$
This allows us to represent $J(x;u(\cd))$ as a quadratic functional on the Hilbert space $L_\dbF^2(\dbR^m)$.

\bp{rep-cost}\sl Suppose that $[A,C]$ is $L^2$-stable. Then there exist a bounded self-adjoint
linear operator $M_2:L_\dbF^2(\dbR^m)\to L_\dbF^2(\dbR^m)$, a bounded linear operator
$M_1:\dbR^n\to L_\dbF^2(\dbR^m)$, an $M_0\in\dbS^n$, and $\h u(\cd)\in L_\dbF^2(\dbR^m)$,
$\h x\in\dbR^n$, $c\in\dbR$ such that
\bel{J-rep}\ba{rll}
\ds J(x;u(\cd))\4n&=\4n&\ds\lan M_2u,u\ran+2\lan M_1x,u\ran+\lan M_0x,x\ran
+2\lan u,\h u\ran+2\lan x,\h x\ran+c,\\
\ns\ds J^0(x;u(\cd))\4n&=\4n&\ds\lan M_2u,u\ran+2\lan M_1x,u\ran+\lan M_0x,x\ran,\\
\ns\4n&~\4n&\ds\qq\q\1n~\forall(x,u(\cd))\in\dbR^n\times L_\dbF^2(\dbR^m). \ea\ee
\ep

\begin{proof}
Let us denote by $X_x^0(\cd)$ the solution of
$$\left\{\2n\ba{ll}
\ds dX(t)=AX(t)dt+CX(t)dW(t),\qq t\ges0, \\
\ns\ds X(0)=x,\ea\right.$$
and by $X_0^u(\cd)$ the solution of
$$\left\{\2n\ba{ll}
\ds dX(t)=\big[AX(t)+Bu(t)\big]dt+\big[CX(t)+Du(t)\big]dW(t),\qq t\ges0, \\
\ns\ds X(0)=0.\ea\right.$$
Define the following linear operators: For any $x\in\dbR^n$ and any $u(\cd)\in L_\dbF^2(\dbR^m)$,
$$(\G x)(\cd)\deq X_x^0(\cd),\qq (Lu)(\cd)\deq X_0^u(\cd).$$
Since $[A,C]$ is $L^2$-stable, Lemma \ref{COCV1} implies that (noting $\cX[0,\i)\subseteq L_\dbF^2(\dbR^n)$)
$$\G:\dbR^n\to L_\dbF^2(\dbR^n),\qq L:L_\dbF^2(\dbR^m)\to L_\dbF^2(\dbR^n)$$
are bounded operators, and that the solution $h(\cd)$ to the SDE
$$\left\{\2n\ba{ll}
\ds dX(t)=\big[AX(t)+b(t)\big]dt+\big[CX(t)+\si(t)\big]dW(t),\qq t\ges0, \\
\ns\ds X(0)=0,\ea\right.$$
is in $\cX[0,\i)$. Let $\G^*$ and $L^*$ denote the adjoint operators of $\G$ and $L$, respectively.
Then, for any $(x,u(\cd))\in\dbR^n\times L_\dbF^2(\dbR^m)$, the corresponding state process $X(\cd)$ is given by
$$X(\cd)=(\G x)(\cd)+(Lu)(\cd)+h(\cd),$$
and hence the cost functional can be written as follows:
\begin{eqnarray*}
J(x;u(\cd)) &\4n=&\4n \lan Q(\G x+Lu+h),\G x+Lu+h\ran+2\lan S(\G x+Lu+h),u\ran\\
&\4n~&\4n +\,\lan Ru,u\ran+2\lan q,\G x+Lu+h\ran+2\lan\rho,u\ran\\
&\4n=&\4n \blan\big(L^*QL+SL+L^*S^\top+R\big)u,u\bran+2\lan(L^*Q+S)\G x,u\ran+\lan\G^*Q\G x,x\ran\\
&\4n~&\4n +\,2\lan L^*(Qh+q)+Sh+\rho,u\ran+2\lan\G^*(Qh+q),x\ran+\lan Qh+2q,h\ran\\
&\4n\equiv&\4n \lan M_2u,u\ran+2\lan M_1x,u\ran+\lan M_0x,x\ran+2\lan u,\h u\ran+2\lan x,\h x\ran+c.
\end{eqnarray*}
In the above, $\lan\cd\,,\cd\ran$ is used for inner products in possibly different spaces.
Note that $\h u(\cd)$, $\h x$ and $c$ vanish in the case when $b(\cd)$, $\si(\cd)$, $q(\cd)$, $\rho(\cd)=0$.
This gives \rf{J-rep}.
\end{proof}

The representation \rf{J-rep} of the cost functional has several consequences,
which we summarize as follows.

\bp{prop123}\sl Suppose that $[A,C]$ is $L^2$-stable. We have the following results:\\[-1.6em]
\begin{enumerate}[\indent\rm(i)]
\item Problem {\rm(LQ)} is open-loop solvable at $x$ if and only if $M_2\ges0$ and
$M_1x+\h u\in\sR( M_2)$. In this case, $u^*(\cd)$ is an open-loop optimal control for
the initial state $x$ if and only if $M_2u^*+M_1x+\h u=0$. \\[-1.6em]
\item If Problem {\rm(LQ)} is open-loop solvable, then so is Problem {\rm(LQ)$^0$}. \\[-1.6em]
\item If Problem {\rm(LQ)$^0$} is open-loop solvable, then there exists a
$U^*(\cd)\in L_\dbF^2(\dbR^{m\times n})$ such that for any $x\in\dbR^n$,
$U^*(\cd)x$ is an open-loop optimal control for the initial state $x$.
\end{enumerate}\ep

\begin{proof}
(i) By definition, a $u^*(\cd)\in L_\dbF^2(\dbR^m)$ is an open-loop
optimal control for the initial state $x$ if and only if
\bel{10.25-1}J(x;u^*(\cd)+\l v(\cd))-J(x;u^*(\cd))\ges0,
\qq\forall v(\cd)\in L_\dbF^2(\dbR^m),~\forall \l\in\dbR^n.\ee
From \rf{J-rep} we have
\begin{eqnarray*}
&&J(x;u^*(\cd)+\l v(\cd))\\
&&=\lan M_2(u^*+\l v),u^*+\l v\ran+2\lan M_1x,u^*+\l v\ran+\lan M_0x,x\ran
+2\lan u^*+\l v,\h u\ran+2\lan x,\h x\ran+c\\
&&=\lan M_2u^*,u^*\ran+2\l\lan M_2u^*,v\ran+\l^2\lan M_2v,v\ran+2\lan M_1x,u^*\ran
+2\l\lan M_1x,v\ran+\lan M_0x,x\ran\\
&&~~~+2\lan u^*,\h u\ran+2\l\lan v,\h u\ran+2\lan x,\h x\ran+c\\
&&=J(x;u^*(\cd))+\l^2\lan M_2v,v\ran+2\l\lan M_2u^*+M_1x+\h u,v\ran.
\end{eqnarray*}
Thus, \rf{10.25-1} is equivalent to
$$\l^2\lan M_2v,v\ran+2\l\lan M_2u^*+M_1x+\h u,v\ran\ges0,
\qq\forall v(\cd)\in L_\dbF^2(\dbR^m),~\forall \l\in\dbR^n,$$
which in turn is equivalent to
$$\lan M_2v,v\ran\ges0,\q\forall v(\cd)\in L_\dbF^2(\dbR^m)\qq \hb{and}\qq M_2u^*+M_1x+\h u=0.$$
The conclusions follow readily.

\ms

(ii) If Problem {\rm(LQ)} is open-loop solvable, then we have by (i): $M_2\ges0$, and
$$M_1x+\h u\in\sR( M_2),\qq\forall x\in\dbR^n.$$
In particular, by taking $x=0$, we see that $\h u\in\sR( M_2)$,
and hence $M_1x\in\sR( M_2)$ for all $x\in\dbR^n$. Using (i) again,
we obtain the open-loop solvability of Problem {\rm(LQ)$^0$}.

\ms

(iii) Let $e_1,\cdots,e_n$ be the standard basis for $\dbR^n$,
and let $u^*_i(\cd)$ be an open-loop optimal control for the initial state $e_i$.
Then $U^*(\cd)\deq(u^*_1(\cd),\cdots,u^*_n(\cd))$ has the desired properties.
\end{proof}

Let us observe that $M_2\ges0$ if and only if
\bel{convex}J^0(0;u(\cd))\ges0,\qq\forall u(\cd)\in L_\dbF^2(\dbR^m).\ee
Further, if there exists a constant $\d>0$ such that $M_2\ges\d I$, or equivalently,
\bel{uni-con}J^0(0;u(\cd))\ges\d\,\dbE\int_0^\i|u(t)|^2dt,\qq\forall u(\cd)\in L_\dbF^2(\dbR^m),\ee
then, by Proposition \ref{prop123} (i), Problem (LQ)$^0$ is uniquely solvable,
with the unique optimal control for the initial state $x$ given by
$$u_x^*(\cd)=-M_2^{-1}M_1x.$$
Note that the value function of Problem (LQ)$^0$ is now given by
\bel{10.28-V0}V^0(x)=\blan\big(M_0-M_1^*M_2^{-1}M_1\big)x,x\bran\equiv\lan Px,x\ran,\qq x\in\dbR^n.\ee
The following result shows that $P$ defined in \rf{10.28-V0} is a static stabilizing solution of ARE \rf{ARE}.

\bt{bt-uni-tu}\sl Suppose that $[A,C]$ is $L^2$-stable and that \rf{uni-con} holds for some $\d>0$.
Then the matrix $P$ defined in \rf{10.28-V0} solves the ARE
\bel{ARE-S}\left\{\2n\ba{ll}
\ds PA+A^\top P+C^\top PC+Q\\
\ns\ds\hphantom{PA}
-\big(PB+C^\top PD+S^\top\big)\big(R+D^\top PD\big)^{-1}\big(B^\top P+D^\top PC+S\big)=0,\\
\ns\ds R+D^\top PD>0,\ea\right.\ee
and
\bel{9-20-Th}\Th\deq-\big(R+D^\top PD\big)^{-1}\big(B^\top P+D^\top PC+S\big)\ee
is a stabilizer of $[A,C;B,D]$. Moreover, the unique open-loop optimal control of Problem {\rm(LQ)$^0$}
for the initial state $x$ is given by
$$u_x^*(\cd)=\Th X_\Th(\cd\,;x),$$
where $X_\Th(\cd\,;x)$ is the solution to the following closed-loop system:
$$\left\{\2n\ba{ll}
\ds dX(t)=(A+B\Th)X(t)dt+(C+D\Th)X(t)dW(t),\qq t\ges0, \\
\ns\ds X(0)= x.\ea\right.$$
\et

\begin{proof}
Let $\F(\cd)$ be the solution of \rf{F}.
Since $[A,C]$ is $L^2$-stable, the following is well-defined:
$$G\deq\dbE\int_0^\i\F(t)^\top Q\F(t)dt$$
For $T>0$, let us consider the state equation
\bel{state-T}\left\{\2n\ba{ll}
\ds dX_T(t)=\big[AX_T(t)+Bu(t)\big]dt+\big[CX_T(t)+Du(t)\big]dW(t),\qq t\in[0,T],\\
\ns\ds X_T(0)=x,\ea\right.\ee
and the cost functional
$$ J_T(x;u(\cd))\deq\dbE\left\{\lan GX_T(T),X_T(T)\ran
+\int_0^T\llan\begin{pmatrix}Q&S^\top\\S&R\\\end{pmatrix}
              \begin{pmatrix}X_T(t)\\ u(t)\end{pmatrix},
              \begin{pmatrix}X_T(t)\\ u(t)\end{pmatrix}\rran dt\right\}.$$
We claim that
\bel{JT-con}J_T(0;u(\cd))\ges\d\,\dbE\int_0^T|u(t)|^2dt,\qq\forall u(\cd)\in L_\dbF^2(0,T;\dbR^m).\ee
To show this, take any $u(\cd)\in L_\dbF^2(0,T;\dbR^m)$ and let $X_T(\cd)$ be the corresponding
solution to \rf{state-T} with initial state $x$. Define the {\it zero-extension} of $u(\cd)$ as follows:
$$[u(\cd)\oplus0{\bf1}_{(T,\i)}](s)=\left\{\2n\ba{ll}
\ds u(t),\q&\ds t\in[0,T],\\
\ns\ds 0,\q&\ds t\in(T,\i).
\ea\right.$$
Then $v(\cd)\equiv[u(\cd)\oplus0{\bf1}_{(T,\i)}](\cd)\in L_\dbF^2(\dbR^m)$, and the solution $X(\cd)$ of
$$\left\{\2n\ba{ll}
\ds dX(t)=\big[AX(t)+Bv(t)\big]dt+\big[CX(t)+Dv(t)\big]dW(t),\qq t\ges0,\\
\ns\ds X(0)=x,
\ea\right.$$
satisfies
$$X(t)=\left\{\2n\ba{ll}
\ds X_T(t),                  \q&\ds t\in[0,T],\\
\ns\ds \F(t)\F(T)^{-1}X_T(T),\q&\ds t\in(T,\i).
\ea\right.$$
Note that for $t\ges T$, $\F(t)\F(T)^{-1}$ has the same distribution as $\F(t-T)$
and is independent of $\cF_T$. Thus,
\bel{9-19-JT}\ba{lll}
\ds J_T(x;u(\cd))\4n&=\4n&\ds
\dbE\bigg\{\bigg\langle\bigg(\dbE\int_0^\i\F(t)^\top Q\F(t)dt\bigg)X_T(T),X_T(T)\bigg\rangle\\
\ns\ds\4n&~\4n&\qq\ds
+\int_0^T\llan\begin{pmatrix}Q&S^\top\\ S&R\end{pmatrix}
              \begin{pmatrix}X_T(t)\\ u(t)\end{pmatrix},
              \begin{pmatrix}X_T(t)\\ u(t)\end{pmatrix}\rran dt\bigg\}\\
\ns\ds\4n&=\4n&\ds
\dbE\bigg\{\bigg\langle\bigg(\dbE\int_T^\i\[\F(t)\F(T)^{-1}\]^\top
Q\[\F(t)\F(T)^{-1}\]dt\bigg) X_T(T),X_T(T)\bigg\rangle\\
\ns\4n&~\4n&\qq\ds
+\int_0^T\llan\begin{pmatrix}Q&S^\top\\ S&R\\\end{pmatrix}
              \begin{pmatrix}X_T(t)\\ u(t)\end{pmatrix},
              \begin{pmatrix}X_T(t)\\ u(t)\end{pmatrix}\rran dt\bigg\}\\
\ns\4n&=\4n&\ds
\dbE\bigg\{\int_T^\i\blan Q\F(t)\F(T)^{-1}X_T(T),\F(t)\F(T)^{-1}X_T(T)\bran dt\\
\ns\4n&~\4n&\qq\ds
+\int_0^T\llan\begin{pmatrix}Q&S^\top\\S&R\\ \end{pmatrix}
              \begin{pmatrix}X_T(t)\\ u(t)\end{pmatrix},
              \begin{pmatrix}X_T(t)\\ u(t)\end{pmatrix}\rran dt\bigg\}\\
\ns\4n&=\4n&\ds \dbE\bigg\{\int_T^\i\lan QX(t),X(t)\ran dt
+\int_0^T\llan\begin{pmatrix}Q&S^\top\\S&R\end{pmatrix}
              \begin{pmatrix}X(t)\\ v(t)\end{pmatrix},
              \begin{pmatrix}X(t)\\ v(t)\end{pmatrix}\rran dt\bigg\}\\
\ns\4n&=\4n&\ds J^0(x;[u(\cd)\oplus0{\bf1}_{(T,\i)}](\cd)).
\ea\ee
In particular, taking $x=0$, we obtain
$$ J_T(0;u(\cd))=J^0(0;[u(\cd)\oplus0{\bf1}_{(T,\i)}](\cd))\ges
\d\,\dbE\int_0^\i\big|[u(\cd)\oplus0{\bf1}_{(T,\i)}](t)\big|^2dt=\d\,\dbE\int_0^T|u(t)|^2dt.$$
This proves our claim.

\ms

The fact \rf{JT-con} allows us to use \cite[Theorem 4.6]{Sun-Li-Yong 2016} to conclude that
for any $T>0$, the differential Riccati equation
$$\left\{\2n\ba{ll}
\ds\dot P(t;T)+P(t;T)A+A^\top P(t;T)+C^\top P(t;T)C+Q\\
\ns\ds\q-\lt[P(t;T)B+C^\top P(t;T)D+S^\top\rt]\lt[R+D^\top P(t;T)D\rt]^{-1}\\
\ns\ds\qq~\cd\lt[B^\top P(t;T)+D^\top P(t;T)C+S\rt]=0,\qq t\in[0,T],\\
\ns\ds P(T;T)=G
\ea\right.$$
admits a unique solution $P(\cd\,;T)\in C([0,T];\dbS^n)$ such that
\begin{eqnarray}
&\label{R+DPD-T}
R+D^\top P(t;T)D\ges\d I,\qq\forall t\in[0,T],&\\
&\label{VT}\ds
\inf_{u(\cd)\in L_\dbF^2(0,T;\dbR^m)}J_T(x;u(\cd))=\lan P(0;T)x,x\ran,\qq\forall x\in\dbR^n.&
\end{eqnarray}
We are going to show that $\{P(0;T)\}_{T>0}$ converges to $P$ as $T\to\i$.
To this end, one observes that \rf{9-19-JT} implies
$$\lan Px,x\ran\les J^0(x;[u(\cd)\oplus0{\bf1}_{(T,\i)}](\cd))=J_T(x;u(\cd)),
\qq\forall x\in\dbR^n,~\forall u(\cd)\in L_\dbF^2(0,T;\dbR^m).$$
Taking infimum over $u(\cd)\in L_\dbF^2(0,T;\dbR^m)$, we obtain
\bel{9-20-0}\lan Px,x\ran\les\lan P(0;T)x,x\ran,\qq\forall x\in\dbR^n,~\forall\, T>0.\ee
On the other hand, for any given $\e>0$, one can find a $u^\e(\cd)\in L_\dbF^2(\dbR^m)$ such that
\bel{9-20-1}\dbE\int_0^\i\llan\begin{pmatrix}Q&S^\top\\S&R\\\end{pmatrix}
                              \begin{pmatrix}X^\e(t)\\ u^\e(t)\end{pmatrix},
                              \begin{pmatrix}X^\e(t)\\ u^\e(t)\end{pmatrix}\rran dt
=J^0(x;u^\e(\cd))\les \lan Px,x\ran+\e,\ee
where $X^\e(\cd)$ is the solution of
$$\left\{\2n\ba{ll}
\ds dX^\e(t)=\big[AX^\e(t)+Bu^\e(t)\big]dt+\big[CX^\e(t)+Du^\e(t)\big]dW(t),\qq t\ges0, \\
\ns\ds X^\e(0)= x.\ea\right.$$
Since by Lemma \ref{COCV1} $X^\e(\cd)\in\cX[0,\i)$,  we have for large $T>0$,
$$\lt|\dbE\lan GX^\e(T),X^\e(T)\ran\rt|\les\e,\qq
\lt|\dbE\int_T^\i\llan\begin{pmatrix}Q&S^\top\\ S&R\end{pmatrix}
                      \begin{pmatrix}X^\e(t) \\ u^\e(t)\end{pmatrix},
                      \begin{pmatrix}X^\e(t) \\ u^\e(t)\end{pmatrix}\rran dt\rt|\les\e.$$
Take $u^\e_T(\cd)=u^\e(\cd)|_{[0,T]}$. Then
\bel{9-20-2}\ba{lll}
\ds J^0(x;u^\e(\cd))\4n&=\4n&\ds J_T(x;u^\e_T(\cd))-\dbE\lan GX^\e(T),X^\e(T)\ran\\
\ns\4n&~\4n&\ds+\,\dbE\int_T^\i\llan\begin{pmatrix}Q&S^\top\\S&R\\\end{pmatrix}
                    \begin{pmatrix}X^\e(t)\\ u^\e(t)\end{pmatrix},
                    \begin{pmatrix}X^\e(t)\\ u^\e(t)\end{pmatrix}\rran dt\ges\lan P(0;T)x,x\ran-2\e.
\ea\ee
Combining \rf{9-20-1} and \rf{9-20-2}, we see that for large $T>0$,
$$\lan P(0;T)x,x\ran\les \lan Px,x\ran+3\e,$$
which, together with \rf{9-20-0}, implies that $P(0;T)\to P$ as $T\to\i$.
Now it follows from \rf{R+DPD-T} that $R+D^\top PD>0$.
Thus, by Lemma \ref{bl-useful lemma}, $P$ solves ARE \rf{ARE-S}.

\ms

Finally, let $\Th$ be as in \rf{9-20-Th}, and for any initial state $x$,
let $(X^*_x(\cd),u_x^*(\cd))$ be the corresponding optimal pair of Problem (LQ)$^0$.
By applying It\^o's formula to $t\to\lan PX^*_x(t),X^*_x(t)\ran$, we have
\begin{eqnarray*}
\lan Px,x\ran\4n&=\4n& J^0(x,u_x^*(\cd))\\
\4n&=\4n& \lan Px,x\ran+\dbE\int_0^\i\[\blan\big(PA+A^\top P+C^\top PC+Q\big)X^*_x(t),X^*_x(t)\bran\\
\4n&~\4n&\hphantom{\lan Px,x\ran+\dbE\int_0^\i\[}
+2\blan\big(B^\top\1n P\1n+\1n D^\top\1n PC\1n+\1n S\big)X^*_x(t),u_x^*(t)\bran\1n
+\1n\blan\big(R\1n+\1n D^\top\1n PD\big)u_x^*(t),u_x^*(t)\bran\]dt\\
\4n&=\4n& \lan Px,x\ran+\dbE\int_0^\i\blan\big(R+D^\top PD\big)
\big[u_x^*(t)-\Th X^*_x(t)\big],u_x^*(t)-\Th X^*_x(t)\bran dt.
\end{eqnarray*}
Since $R+D^\top PD>0$, we must have
$$u_x^*(\cd)=\Th X^*_x(\cd),\qq\forall x\in\dbR^n,$$
and hence $X^*_x(\cd)$ satisfies
$$\left\{\2n\ba{ll}
\ds dX^*_x(t)=(A+B\Th)X^*_x(t)dt+(C+D\Th)X^*_x(t)dW(t),\qq t\ges0, \\
\ns\ds X^*_x(0)= x.\ea\right.$$
Since $X^*_x(\cd)\in\cX[0,\i)$ for all $x\in\dbR^n$,
we conclude that $\Th$ is a stabilizer of $[A,C;B,D]$. The rest of the proof is clear.
\end{proof}

We are now ready to prove Theorem \ref{main} for the case when $[A,C]$ is $L^2$-stable.

\begin{proof}[\indent\textbf{Proof of Theorem {\rm\ref{main}} under the assumption that $[A,C]$ is $L^2$-stable}]
(ii) $\Ra$ (i) is obvious.

\ms

(i) $\Ra$ (iii): First, by Proposition \ref{prop123}, $M_2\ges0$ and Problem (LQ)$^0$ is
also open-loop solvable. For any $\e>0$, let us consider the state equation
$$\left\{\2n\ba{ll}
\ds dX(t)=\big[AX(t)+Bu(t)\big]dt+\big[CX(t)+Du(t)\big]dW(t),\qq t\ges0, \\
\ns\ds X(0)= x,\ea\right.$$
and the cost functional
$$J^0_\e(x; u(\cd)) = J^0(x; u(\cd))+\e\dbE\int_0^\i|u(t)|^2dt
                    = \dbE\int_0^\i\llan\begin{pmatrix}Q&S^\top\\S&R+\e I\end{pmatrix}
                                        \begin{pmatrix}X(t)\\ u(t)\end{pmatrix},
                                        \begin{pmatrix}X(t)\\ u(t)\end{pmatrix}\rran dt.$$
Denote by Problem (LQ)$_\e^0$ the above problem,
and by $V_\e^0(\cd)$ the corresponding value function. Since
$$J^0_\e(0; u(\cd))=\lan(M_2+\e I)u,u\ran\ges\e\dbE\int_0^\i|u(t)|^2dt,
\qq\forall u(\cd)\in L_\dbF^2(\dbR^m),$$
by Theorem \ref{bt-uni-tu}, the following ARE
\bel{ARE-S-e}\left\{\2n\ba{ll}
\ds P_\e A+A^\top P_\e+C^\top P_\e C+Q\\
\ns\ds~-\big(P_\e B+C^\top P_\e D+S^\top\big)\big(R+\e I+D^\top P_\e D\big)^{-1}
\big(B^\top P_\e+D^\top P_\e C+S\big)=0,\\
\ns\ds R+\e I+D^\top P_\e D>0,
\ea\right.\ee
admits a unique solution $P_\e\in\dbS^n$ such that
$V_\e^0(x)=\lan P_\e x,x\ran$ for all $x\in\dbR^n$. Moreover,
$$\Th_\e\deq-\big(R+\e I+D^\top P_\e D\big)^{-1}\big(B^\top P_\e+D^\top P_\e C+S\big)$$
is a stabilizer of $[A,C;B,D]$, and the unique open-loop optimal control $u^*_\e(\cd\,;x)$
of Problem (LQ)$_\e^0$ for the initial state $x$ is given by
$$ u^*_\e(t;x)=\Th_\e\Psi_\e(t)x,\qq t\ges0, $$
where $\Psi_\e(\cd)$ is the solution to the following SDE for $\dbR^{n\times n}$-valued processes:
$$\left\{\2n\ba{ll}
\ds d\Psi_\e(t)=(A+B\Th_\e)\Psi_\e(t)dt+(C+D\Th_\e)\Psi_\e(t)dW(t),\qq t\ges0, \\
\ns\ds \Psi_\e(0)=I.\ea\right.$$
Now let $U^*(\cd)\in L_\dbF^2(\dbR^{m\times n})$ be a process with the property in
Proposition \ref{prop123} (iii). By the definition of value function, we have
for any $x\in\dbR^n$ and $\e>0$,
\bel{10-14-2}\ba{lll}
\ds V^0(x)+\e\dbE\int_0^\i|\Th_\e\Psi_\e(t)x|^2dt
\4n&\les\4n&\ds J^0(x;\Th_\e\Psi_\e(\cd)x)+\e\dbE\int_0^\i|\Th_\e\Psi_\e(t)x|^2dt\\
\ns\4n&=\4n&\ds J_\e^0(x;\Th_\e\Psi_\e(\cd)x)=V^0_\e(x)\les J_\e^0(x;U^*(\cd)x)\\
\ns\4n&=\4n&\ds V^0(x)+\e\dbE\int_0^\i|U^*(t)x|^2dt,
\ea\ee
which implies the following:
\begin{eqnarray}
&\label{10-13-1}\ds
V^0(x)\les V^0_\e(x)=\lan P_\e x,x\ran\les V^0(x)+\e\dbE\int_0^\i|U^*(t)x|^2dt,
\qq\forall x\in\dbR^n,~\forall\e>0,&\\
&\label{10-13-2}\ds
0\les\dbE\int_0^\i\Psi_\e(t)^\top\Th_\e^\top\Th_\e\Psi_\e(t)dt
\les\dbE\int_0^\i U^*(t)^\top U^*(t)dt,\qq\forall\e>0.&
\end{eqnarray}
From \rf{10-13-1} we see that $P\equiv\ds\lim_{\e\to0}P_\e$ exists and $V^0(x)=\lan Px,x\ran$
for all $x\in\dbR^n$. Denote
$$\Pi_\e=\dbE\int_0^\i\Psi_\e(t)^\top\Th_\e^\top\Th_\e\Psi_\e(t)dt.$$
It follows from \rf{10-13-2} that $\{\Pi_\e\}_{\e>0}$ is bounded.
Moreover, noting that $[A+B\Th_\e,C+D\Th_\e] $ is $L^2$-stable,
we have by Lemma \ref{AC-stable} that
$$\Pi_\e(A+B\Th_\e)+(A+B\Th_\e)^\top\Pi_\e+(C+D\Th_\e)^\top\Pi_\e(C+D\Th_\e)+\Th_\e^\top\Th_\e=0.$$
Thus,
\begin{eqnarray*}
0\les\Th_\e^\top\Th_\e
\4n&=\4n& -\,\big[\Pi_\e(A+B\Th_\e)+(A+B\Th_\e)^\top\Pi_\e+(C+D\Th_\e)^\top\Pi_\e(C+D\Th_\e)\big]\\
\4n&\les\4n& -\,\big[\Pi_\e(A+B\Th_\e)+(A+B\Th_\e)^\top\Pi_\e\big],\qq\forall\e>0.
\end{eqnarray*}
The above, together with the boundedness of $\{\Pi_\e\}_{\e>0}$, shows that
\bel{bound-Th-e}|\Th_\e|^2\les K(1+|\Th_\e|),\qq\forall\e>0,\ee
for some constant $K>0$. Noting that \rf{bound-Th-e} implies the boundedness of $\{\Th_\e\}_{\e>0}$,
we may choose a sequence $\{\e_k\}^\i_{k=1}\subseteq(0,\i)$ with $\ds\lim_{k\to\i}\e_k=0$ such that
$\Th\equiv\ds\lim_{k\to\i}\Th_{\e_k}$ exists. Then
\bel{10.15-1}\big(R+D^\top P D\big)\Th=\lim_{k\to\i}\big(R+\e_k I+D^\top P_{\e_k} D\big)\Th_{\e_k}
=-\big(B^\top P+D^\top PC+S\big).\ee
It follows from Lemma \ref{Penrose} that
\bel{10-14-1}\left\{\2n\ba{ll}
\ds \sR\big(B^\top P+D^\top PC+S\big)\subseteq\sR\big(R+D^\top PD\big),\\
\ns\ds \Th=-\cN(P)^\dag\cL(P)^\top+\big[I-\cN(P)^\dag\cN(P)\big]\Pi\q\hb{for some}~\Pi\in\dbR^{m\times n}.
\ea\right.\ee
A passage to the limit along $\{\e_k\}^\i_{k=1}$ in \rf{ARE-S-e} yields
\bel{10.15-2}\left\{\2n\ba{ll}
\ds PA+A^\top P+C^\top PC+Q-\Th^\top\big(R+D^\top PD\big)\Th=0,\\
\ns\ds R+D^\top PD\ges0,
\ea\right.\ee
which, together with \rf{10-14-1}, implies that $P$ solves ARE \rf{ARE}.
To see that $P$ is a static stabilizing solution, we need only show that $\Th\in\sS[A,C;B,D]$.
For this, let $\Psi(\cd)$ be the solution of
$$\left\{\2n\ba{ll}
\ds d\Psi(t)=(A+B\Th)\Psi(t)dt+(C+D\Th)\Psi(t)dW(t),\qq t\ges0, \\
\ns\ds \Psi(0)=I.\ea\right.$$
Since $\Th_{\e_k}\to\Th$ as $k\to\i$, we have $\Psi_{\e_k}(t)\to\Psi(t),~\as$ for all $t\ges0$.
Using Fatou's lemma and \rf{10-14-2}, we obtain
$$\dbE\int_0^\i|\Th\Psi(t)x|^2dt\les\liminf_{k\to\i}\dbE\int_0^\i|\Th_{\e_k}\Psi_{\e_k}(t)x|^2dt
\les\dbE\int_0^\i|U^*(t)x|^2dt<\i,\qq\forall x\in\dbR^n.$$
This implies $\Th\Psi(\cd)\in L_\dbF^2(\dbR^{m\times n})$, and thus by Lemma \ref{COCV1}
$\Psi(\cd)\in L_\dbF^2(\dbR^{n\times n})$. Therefore, $\Th$ is a stabilizer of $[A,C;B,D]$.

\ms

Next, consider the following BSDE on $[0,\i)$:
\bel{10-14-3}d\eta(t)=-\big[(A+B\Th)^\top\eta+(C+D\Th)^\top\z+(C+D\Th)^\top P\si
+\Th^\top\rho+Pb+q\big]dt+\z dW(t).\ee
Since $[A+B\Th,C+D\Th]$ is $L^2$-stable, it follows form Lemma \ref{COCV2} that \rf{10-14-3}
admits a unique $L^2$-stable adapted solution $(\eta(\cd),\z(\cd))$. For any initial state
$x$ and any control process $u(\cd)\in L_\dbF^2(\dbR^m)$, let $X(\cd)\equiv X(\cd\,;x,u(\cd))$
be the corresponding solution of \rf{state}.
Applying It\^o's formula to $t\mapsto\lan PX(t),X(t)\ran$, we have
\bel{Ito-1}\ba{lll}
\ds-\lan Px,x\ran\4n&=\4n&\ds \dbE\int_0^\i\[\lan P(AX+Bu+b),X\ran+\lan PX,AX+Bu+b\ran\\
\ns\4n&~\4n&\ds\hphantom{\dbE\int_0^\i\[}
+\lan P(CX+Du+\si),CX+Du+\si\ran\]dt\\
\ns\4n&=\4n&\ds\dbE\int_0^\i\[\blan\big(PA+A^\top P+C^\top PC\big)X,X\bran
+2\blan\big(B^\top P+D^\top PC\big)X,u\bran\\
\ns\4n&~\4n&\ds\hphantom{\dbE\int_0^\i\[}
+\blan D^\top PDu,u\bran+2\blan C^\top P\si+Pb,X\bran
+2\blan D^\top P\si,u\bran+\lan P\si,\si\ran\]dt.\ea\ee
Applying It\^o's formula to $t\mapsto\lan \eta(t),X(t)\ran$, we have
\bel{Ito-2}\ba{lll}
\ds\dbE\lan\eta(0),x\ran\4n&=\4n&\ds
\dbE\int_0^\i\[\blan(A\1n+\1nB\Th)^\top\eta+(C\1n+\1nD\Th)^\top\z
+(C\1n+\1nD\Th)^\top P\si+\Th^\top\rho+Pb+q,X\bran\\
\ns\4n&~\4n&\ds\hphantom{\dbE\int_0^\i\[}
-\lan\eta,AX+Bu+b\ran-\lan\z,CX+Du+\si\ran\]dt\\
\ns\4n&=\4n&\ds
\dbE\int_0^\i\[\blan\Th^\top\big(B^\top\eta+D^\top\z+D^\top P\si+\rho\big),X\bran
+\blan C^\top P\si+Pb+q,X\bran\\
\ns\4n&~\4n&\ds\hphantom{\dbE\int_0^\i\[}
-\blan B^\top\eta+D^\top\z,u\bran-\lan\eta,b\ran-\lan\z,\si\ran\]dt.\ea\ee
Then using \rf{10.15-1} and \rf{10.15-2}, we have
\bel{10.15-3}\ba{lll}
\ds J(x;u(\cd))-\lan Px,x\ran-2\dbE\lan\eta(0),x\ran\\
\ns\ds=\dbE\int_0^\i\[\blan\big(PA+A^\top P+C^\top PC+Q\big)X,X\bran
+2\blan\big(B^\top P+D^\top PC+S\big)X,u\bran\\
\ns\ds\hphantom{=\dbE\int_0^\i\[}
+\blan\big(R+D^\top PD\big)u,u\bran
+2\blan B^\top\eta+D^\top\z+D^\top P\si+\rho,u-\Th X\bran\\
\ns\ds\hphantom{=\dbE\int_0^\i\[}
+\lan P\si,\si\ran+2\lan\eta,b\ran+2\lan\z,\si\ran\]dt\\
\ns\ds=\dbE\int_0^\i\[\blan\big(R+D^\top PD\big)\Th X,\Th X\bran
-2\blan\big(R+D^\top PD\big)\Th X,u\bran+\blan\big(R+D^\top PD\big)u,u\bran\\
\ns\ds\hphantom{=\dbE\int_0^\i\[}
+2\blan B^\top\eta+D^\top\z+D^\top P\si+\rho,u-\Th X\bran
+\lan P\si,\si\ran+2\lan\eta,b\ran+2\lan\z,\si\ran\]dt\\
\ns\ds=\dbE\int_0^\i\[\blan\big(R+D^\top PD\big)(u-\Th X),(u-\Th X)\bran
+2\blan B^\top\eta+D^\top\z+D^\top P\si+\rho,u-\Th X\bran\\
\ns\ds\hphantom{=\dbE\int_0^\i\[}
+\lan P\si,\si\ran+2\lan\eta,b\ran+2\lan\z,\si\ran\]dt. \ea\ee
Let $u^*(\cd)$ be an open-loop optimal control of Problem {\rm(LQ)} for the initial state $x$,
and denote by $X_\Th(\cd\,;x,v(\cd))$ the solution to the following SDE:
$$\left\{\2n\ba{ll}
\ds dX_\Th(t)=\big[(A+B\Th)X_\Th(t)+Bv(t)+b(t)\big]dt\\
\ns\ds\hphantom{dX_\Th(t)=}
+\big[(C+D\Th)X_\Th(t)+Dv(t)+\si(t)\big]dW(t),\qq t\ges0, \\
\ns\ds X_\Th(0)= x.\ea\right.$$
By Proposition \ref{p-ad},
$$ u^*(\cd)=\Th X_\Th(\cd\,;x,v^*(\cd))+v^*(\cd),$$
for some $v^*(\cd)\in L_\dbF^2(\dbR^m)$, and hence
\bel{10.30-1}\ba{lll}
\ds J(x;\Th X_\Th(\cd\,;x,v^*(\cd))+v^*(\cd))\4n&=\4n&\ds J(x;u^*(\cd))\\
\ns\4n&\les\4n&\ds J(x;\Th X_\Th(\cd\,;x,v(\cd))+v(\cd)),
\qq\forall v(\cd)\in L_\dbF^2(\dbR^m).\ea\ee
Taking $u(\cd)=\Th X_\Th(\cd\,;x,v(\cd))+v(\cd)$ and noting that
$$X(\cd\,;x,u^*(\cd))=X_\Th(\cd\,;x,v^*(\cd)),\qq X(\cd\,;x,u(\cd))=X_\Th(\cd\,;x,v(\cd)),$$
we have from \rf{10.15-3} and \rf{10.30-1} that for any $v(\cd)\in L_\dbF^2(\dbR^m)$,
\bel{10.15-4}\ba{lll}
\ds\lan Px,x\ran+2\dbE\lan\eta(0),x\ran
+\dbE\int_0^\i\[\lan P\si,\si\ran+2\lan\eta,b\ran+2\lan\z,\si\ran\]dt\\
\ns\ds\hphantom{=}
+\dbE\int_0^\i\[\blan\big(R+D^\top PD\big)v^*,v^*\bran
+2\blan B^\top\eta+D^\top\z+D^\top P\si+\rho,v^*\bran\]dt\\
\ns\ds=J(x;u^*(\cd))\les J(x;u(\cd))\\
\ns\ds\les\lan Px,x\ran+2\dbE\lan\eta(0),x\ran
+\dbE\int_0^\i\[\lan P\si,\si\ran+2\lan\eta,b\ran+2\lan\z,\si\ran\]dt\\
\ns\ds\hphantom{=}
+\dbE\int_0^\i\[\blan\big(R+D^\top PD\big)v,v\bran
+2\blan B^\top\eta+D^\top\z+D^\top P\si+\rho,v\bran\]dt,
\ea\ee
which shows that $v^*(\cd)$ is a minimizer of the functional
$$F(v(\cd))\deq\dbE\int_0^\i\[\blan\big(R+D^\top PD\big)v,v\bran
+2\blan B^\top\eta+D^\top\z+D^\top P\si+\rho,v\bran\]dt,\q~ v(\cd)\in L_\dbF^2(\dbR^m).$$
Therefore, we must have
$$\big(R+D^\top PD\big)v^*+B^\top\eta+D^\top\z+D^\top P\si+\rho=0.$$
It follows from Lemma \ref{Penrose} that
\bel{10.15-5}\left\{\2n\ba{ll}
\ds B^\top\eta+D^\top\z+D^\top P\si+\rho\in\sR\big(R+D^\top PD\big),\\
\ns\ds v^*=-\cN(P)^\dag\big(B^\top\eta+D^\top\z+D^\top P\si+\rho\big)
+\big[I-\cN(P)^\dag\cN(P)\big]\n\q\hb{for some}~\n\in L_\dbF^2(\dbR^m).
\ea\right.\ee
Recall \rf{10-14-1} and observe that
$\big[\Th^\top\1n+\cL(P)\cN(P)^\dag\big]\big(B^\top\eta+D^\top\z+D^\top P\si+\rho\big)=0$. Thus,
\begin{eqnarray*}
&&(A+B\Th)^\top\eta+(C+D\Th)^\top\z+(C+D\Th)^\top P\si+\Th^\top\rho+Pb+q\\
&&=A^\top\eta+C^\top\z+C^\top P\si+Pb+q
+\Th^\top\big(B^\top\eta+D^\top\z+D^\top P\si+\rho\big)\\
&&=A^\top\eta+C^\top\z+C^\top P\si+Pb+q
-\cL(P)\cN(P)^\dag\big(B^\top\eta+D^\top\z+D^\top P\si+\rho\big)\\
&&=\big[A^\top-\cL(P)\cN(P)^\dag B^\top\big]\eta+\big[C^\top-\cL(P)\cN(P)^\dag D^\top\big]\z\\
&&\hphantom{=}+\big[C^\top-\cL(P)\cN(P)^\dag D^\top\big]P\si-\cL(P)\cN(P)^\dag\rho+Pb+q.
\end{eqnarray*}
We see then $(\eta(\cd),\z(\cd))$ is an $L^2$-stable adapted solution of \rf{eta-zeta}.
Further, combining \rf{10.15-4} and \rf{10.15-5}, we have
\begin{eqnarray*}
V(x)\4n&=\4n& J(x;u^*(\cd))\\
\4n&=\4n& \lan Px,x\ran+2\dbE\lan\eta(0),x\ran
+\dbE\int_0^\i\[\lan P\si,\si\ran+2\lan\eta,b\ran+2\lan\z,\si\ran\]dt\\
\4n&~\4n& +\,\dbE\int_0^\i\[\blan\big(R+D^\top PD\big)v^*,v^*\bran
+2\blan B^\top\eta+D^\top\z+D^\top P\si+\rho,v^*\bran\]dt\\
\4n&=\4n& \lan Px,x\ran+2\dbE\lan\eta(0),x\ran
+\dbE\int_0^\i\[\lan P\si,\si\ran+2\lan\eta,b\ran+2\lan\z,\si\ran\]dt\\
\4n&~\4n& -\,\dbE\int_0^\i\blan \cN(P)^\dag\big(B^\top\eta+D^\top\z+D^\top P\si+\rho\big),
B^\top\eta+D^\top\z+D^\top P\si+\rho\bran dt.
\end{eqnarray*}

\ss

(iii) $\Ra$ (ii): We take any $(x,u(\cd))\in\dbR^n\times L_\dbF^2(\dbR^m)$,
and let $X(\cd)\equiv X(\cd\,;x,u(\cd))$ be the corresponding state process.
Proceeding by analogy with \rf{Ito-1}--\rf{10.15-3}, we obtain
\begin{eqnarray*}
J(x;u(\cd)) \4n&=\4n& \lan Px,x\ran+2\dbE\lan\eta(0),x\ran
+\dbE\int_0^\i\[\lan P\si,\si\ran+2\lan\eta,b\ran+2\lan\z,\si\ran\]dt\\
\4n&~\4n& +\,\dbE\int_0^\i\[\blan\big(PA+A^\top P+C^\top PC+Q\big)X,X\bran\\
\4n&~\4n&\hphantom{+\,\dbE\int_0^\i\[}
+2\blan\big(B^\top P+D^\top PC+S\big)X,u\bran
+\blan\big(R+D^\top PD\big)u,u\bran\\
\4n&~\4n&\hphantom{+\,\dbE\int_0^\i\[}
+2\blan B^\top\eta+D^\top\z+D^\top P\si+\rho,u+\cN(P)^\dag\cL(P)^\top X\bran\]dt.
\end{eqnarray*}
Let $(\Th^*,v^*(\cd))$ be defined by \rf{close-rep}. We have
$$\lt\{\2n\ba{ll}
\ds B^\top P+D^\top PC+S=-(R+D^\top PD)\Th^*=-\cN(P)\Th^*,\\
\ns\ds B^\top\eta+D^\top\z+D^\top P\si+\rho=-(R+D^\top PD)v^*=-\cN(P)v^*,\\
\ns\ds\cN(P)\cN(P)^\dag\cL(P)^\top=-\cN(P)\Th^*.
\ea\rt.$$
Thus,
\bel{11.4-1}\ba{lll}
\ds J(x;u(\cd))\4n&=\4n&\ds\lan Px,x\ran+2\dbE\lan\eta(0),x\ran
+\dbE\int_0^\i\[\lan P\si,\si\ran+2\lan\eta,b\ran+2\lan\z,\si\ran\]dt\\
\ns\4n&~\4n&\ds+\,\dbE\int_0^\i\[\blan(\Th^*)^\top\cN(P)\Th^*X,X\bran
-2\lan\cN(P)\Th^*X,u\ran+\lan\cN(P)u,u\ran\\
\ns\4n&~\4n&\ds\hphantom{+\,\dbE\int_0^\i\[}
-2\blan\cN(P)v^*,u+\cN(P)^\dag\cL(P)^\top X\bran\]dt\\
\ns\4n&=\4n&\ds\lan Px,x\ran+2\dbE\lan\eta(0),x\ran
+\dbE\int_0^\i\[\lan P\si,\si\ran+2\lan\eta,b\ran+2\lan\z,\si\ran\]dt\\
\ns\4n&~\4n&\ds+\,\dbE\int_0^\i\[\lan\cN(P)(u-\Th^*X),u-\Th^*X\ran
-2\lan v^*,\cN(P)u-\cN(P)\Th^*X\ran\]dt\\
\ns\4n&=\4n&\ds\lan Px,x\ran+2\dbE\lan\eta(0),x\ran
+\dbE\int_0^\i\[\lan P\si,\si\ran+2\lan\eta,b\ran+2\lan\z,\si\ran-\lan\cN(P)v^*,v^*\ran\]dt\\
\ns\4n&~\4n&\ds+\,\dbE\int_0^\i\lan\cN(P)(u-\Th^*X-v^*),u-\Th^*X-v^*\ran dt\\
\ns\4n&=\4n&\ds\lan Px,x\ran+2\dbE\lan\eta(0),x\ran
+\dbE\int_0^\i\[\lan P\si,\si\ran+2\lan\eta,b\ran+2\lan\z,\si\ran\\
\ns\4n&~\4n&\ds-\,\blan \cN(P)^\dag\big(B^\top\eta+D^\top\z+D^\top P\si+\rho\big),
B^\top\eta+D^\top\z+D^\top P\si+\rho\bran\]dt\\
\ns\4n&~\4n&\ds+\,\dbE\int_0^\i\lan\cN(P)(u-\Th^*X-v^*),u-\Th^*X-v^*\ran dt.
\ea\ee
Since $\cN(P)\ges0$ and $\Th^*$ is a stabilizer of $[A,C;B,D]$, we see that
\begin{eqnarray*}
J(x;u(\cd))\4n&\ges\4n& \lan Px,x\ran+2\dbE\lan\eta(0),x\ran
+\dbE\int_0^\i\[\lan P\si,\si\ran+2\lan\eta,b\ran+2\lan\z,\si\ran\\
\ns\4n&~\4n&\ds-\,\blan \cN(P)^\dag\big(B^\top\eta+D^\top\z+D^\top P\si+\rho\big),
B^\top\eta+D^\top\z+D^\top P\si+\rho\bran\]dt\\
\4n&=\4n& J(x;\Th^*X^*(\cd)+v^*(\cd)),\qq\forall x\in\dbR^n,~\forall u(\cd)\in L^2_\dbF(\dbR^m).
\end{eqnarray*}
That is, $(\Th^*,v^*(\cd))$ is a closed-loop optimal strategy of Problem (LQ).

\ms

Finally, if $(\Th,v(\cd))$ is a closed-loop optimal strategy, then with $X(\cd)$
denoting the solution of
$$\left\{\2n\ba{ll}
\ds dX(t)=\big[(A+B\Th)X(t)+Bv(t)+b(t)\big]dt\\
\ns\ds\hphantom{dX(t)=}
+\big[(C+D\Th)X(t)+Dv(t)+\si(t)\big]dW(t),\qq t\ges0,\\
\ns\ds X(0)= x,\ea\right.$$
and $u(\cd)\equiv\Th X(\cd)+v(\cd)$ denoting the outcome of $(\Th,v(\cd))$, \rf{11.4-1} implies that
\begin{eqnarray*}
V(x)\4n&=\4n&J(x;u(\cd))\\
\4n&=\4n& \lan Px,x\ran+2\dbE\lan\eta(0),x\ran
+\dbE\int_0^\i\[\lan P\si,\si\ran+2\lan\eta,b\ran+2\lan\z,\si\ran\\
\4n&~\4n&-\,\blan \cN(P)^\dag\big(B^\top\eta+D^\top\z+D^\top P\si+\rho\big),
B^\top\eta+D^\top\z+D^\top P\si+\rho\bran\]dt\\
\4n&~\4n&+\,\dbE\int_0^\i\lan\cN(P)(u-\Th^*X-v^*),u-\Th^*X-v^*\ran dt\\
\4n&=\4n& V(x)+\dbE\int_0^\i\lan\cN(P)(\Th X+v-\Th^*X-v^*),\Th X+v-\Th^*X-v^*\ran dt.
\end{eqnarray*}
Thus,
\begin{eqnarray*}
&&\dbE\int_0^\i\big|\cN(P)^{1\over2}(\Th X+v-\Th^*X-v^*)\big|^2 dt\\
&&=\dbE\int_0^\i\lan\cN(P)(\Th X+v-\Th^*X-v^*),\Th X+v-\Th^*X-v^*\ran dt=0,\qq\forall x\in\dbR^n.
\end{eqnarray*}
It follows that $\cN(P)^{1\over2}(\Th X+v-\Th^*X-v^*)=0$ for all $x\in\dbR^n$, and hence
$$\cN(P)(\Th-\Th^*)X+\cN(P)(v-v^*)=0,\qq\forall x\in\dbR^n.$$
Since the above holds for each $x\in\dbR^n$, and $\Th$, $\Th^*$, $v(\cd)$, and $v^*(\cd)$
are independent of $x$, by subtracting solutions corresponding $x$ and $0$, the latter from
the former, we see that for any $x\in\dbR^n$, the solution $X_0(\cd)$ of
$$\left\{\2n\ba{ll}
\ds dX_0(t)=(A+B\Th)X_0(t)dt+(C+D\Th)X_0(t)dW(t),\qq t\ges0,\\
\ns\ds X_0(0)=x,\ea\right.$$
satisfies $\cN(P)(\Th-\Th^*)X_0=0$, from which we conclude that $\cN(P)(\Th-\Th^*)=0$ and
$\cN(P)(v-v^*)=0$. Hence,
$$\ba{ll}
\ds\cN(P)\Th=\cN(P)\Th^*=-\cL(P)^\top,\\
\ns\ds\cN(P)v=\cN(P)v^*=-\big(B^\top\eta+D^\top\z+D^\top P\si+\rho\big).\ea$$
It follows from Lemma \ref{Penrose} that $(\Th,v(\cd))$ is of the form \rf{close-rep}.
Similarly, if $u^*(\cd)$ is an open-loop optimal control for the initial state $x$,
then with $X^*(\cd)$ denoting the corresponding optimal state process, we have
$$\cN(P)(u^*-\Th^*X^*-v^*)=0,$$
or equivalently,
$$\cN(P)u^*=\cN(P)\Th^*X^*+\cN(P)v^*=-\cL(P)^\top X^*-\big(B^\top\eta+D^\top\z+D^\top P\si+\rho\big).$$
By Lemma \ref{Penrose}, there exists a $\n(\cd)\in L^2_\dbF(\dbR^m)$ such that
\begin{eqnarray*}
u^*\4n&=\4n& -\,\cN(P)^\dag\cL(P)^\top X^*-\cN(P)^\dag\big(B^\top\eta+D^\top\z+D^\top P\si+\rho\big)
+\big[I-\cN(P)^\dag\cN(P)\big]\n\\
\4n&=\4n& \lt\{-\cN(P)^\dag\cL(P)^\top+\big[I-\cN(P)^\dag\cN(P)\big]\Pi\rt\} X^*\\
\4n&~\4n& -\,\cN(P)^\dag\big(B^\top\eta+D^\top\z+D^\top P\si+\rho\big)
+\big[I-\cN(P)^\dag\cN(P)\big](\n-\Pi X^*).
\end{eqnarray*}
In the above, $\Pi\in\dbR^{m\times n}$ is chosen such that
$-\cN(P)^\dag\cL(P)^\top+\big[I-\cN(P)^\dag\cN(P)\big]\Pi\in\sS[A,C;B,D]$.
This shows that $u^*(\cd)$ has the closed-loop representation \rf{u-star}.
\end{proof}

\section{The Proof of Theorem \ref{main}: The General Case}

We turn to the proof of Theorem \ref{main} for the general case $\sS[A,C;B,D]\ne\emptyset$.
The idea is to apply Proposition \ref{p-ad}, thus converting Problem (LQ) into an equivalent
one, in which the corresponding uncontrolled system is $L^2$-stable.

\ms

More precisely, take any $\Si\in\sS[A,C;B,D]$, and consider the state equation
\bel{Th-state}\left\{\2n\ba{ll}
\ds d\wt X(t)=\big[\wt A\wt X(t)+Bv(t)+b(t)\big]dt+\big[\wt C\wt X(t)+Dv(t)+\si(t)\big]dW(t),\qq t\ges0,\\
\ns\ds \wt X(0)= x,\ea\right.\ee
and the cost functional
\bel{Th-cost}\ba{lll}
\ds \wt J(x;v(\cd))\4n&\deq\4n&\ds J(x;\Si\wt X(\cd)+v(\cd))\\
\ns\4n&=\4n&\ds\dbE\int_0^\i
\Bigg[\llan\begin{pmatrix}Q&S^\top\\S&R\\\end{pmatrix}
           \begin{pmatrix}\wt X(t)\\ \Si\wt X(t)+v(t)\end{pmatrix},
           \begin{pmatrix}\wt X(t)\\ \Si\wt X(t)+v(t)\end{pmatrix}\rran\\
\ns\4n&~\4n&\ds\hphantom{\dbE\int_0^\i\Bigg[}
+2\llan\begin{pmatrix}q(t)\\ \rho(t)\end{pmatrix},
       \begin{pmatrix}\wt X(t)\\ \Si\wt X(t)+v(t)\end{pmatrix}\rran\Bigg]dt\\
\ns\4n&=\4n&\ds\dbE\int_0^\i
\Bigg[\llan\begin{pmatrix}\wt Q&\wt S^\top\\\wt S&R\\\end{pmatrix}
           \begin{pmatrix}\wt X(t)\\ v(t)\end{pmatrix},
           \begin{pmatrix}\wt X(t)\\ v(t)\end{pmatrix}\rran
+2\llan\begin{pmatrix}\wt q(t)\\ \rho(t)\end{pmatrix},
       \begin{pmatrix}\wt X(t)\\ v(t)\end{pmatrix}\rran\Bigg]dt,\ea\ee
where
$$\left\{\2n\ba{llr}
\ds    \wt A\4n&=\4n&\ds A+B\Si, \qq \wt C=C+D\Si,\\
\ns\ds \wt Q\4n&=\4n&\ds Q+S^\top\Si+\Si^\top S+\Si^\top R\Si,\\
\ns\ds \wt S\4n&=\4n&\ds S+R\Si, \qq~\wt q=q+\Si^\top\rho.
\ea\right.$$
Note that the system $[\wt A,\wt C]$ is $L^2$-stable. We denote by $\wt X(\cd\,;x,v(\cd))$
the solution of \rf{Th-state} corresponding to $x$ and $v(\cd)$, and by Problem ($\wt{\hb{LQ}}$)
the above problem. The following lists several basic facts about Problem ($\wt{\hb{LQ}}$),
whose proofs are straightforward consequences of Proposition \ref{p-ad}.\\[-1.6em]
\begin{enumerate}[\indent\rm(1)]
\item Problem ($\wt{\hb{LQ}}$) is open-loop solvable at $x\in\dbR^n$ if and only if
      Problem (LQ) is. In this case, $v^*(\cd)$ is an open-loop optimal control of
      Problem ($\wt{\hb{LQ}}$) if and only if $u^*(\cd)\deq v^*(\cd)+\Si\wt X(\cd\,;x,v^*(\cd))$
      is an open-loop optimal control of Problem (LQ).\\[-1.6em]
\item Problem ($\wt{\hb{LQ}}$) is closed-loop solvable if and only if Problem (LQ) is. In this case,
      $(\Si^*,v^*(\cd))$ is a closed-loop optimal strategy of Problem ($\wt{\hb{LQ}}$) if and only
      if $(\Si^*+\Si,v^*(\cd))$ is a closed-loop optimal strategy of Problem (LQ).
\end{enumerate}

\begin{proof}[\indent\textbf{Proof of Theorem {\rm\ref{main}}, general case}]
(ii) trivially implies (i), and the implication (iii) $\Ra$ (ii) can be proved in the same
way as in the case when $[A,C]$ is $L^2$-stable. For the implication (i) $\Ra$ (iii), we
consider Problem ($\wt{\hb{LQ}}$). Since the open-loop solvabilities of Problem (LQ) and
Problem ($\wt{\hb{LQ}}$) are equivalent, and $[\wt A,\wt C]$ is $L^2$-stable, by the result
for the $L^2$-stable case, the following ARE
$$\left\{\2n\ba{ll}
\ds P\wt A+\wt A^\top P+\wt C^\top P\wt C+\wt Q
-\big(PB+\wt C^\top PD+\wt S^\top\big)\big(R+D^\top PD\big)^\dag\big(B^\top P+D^\top P\wt C+\wt S\big)=0,\\
\ns\ds\sR\big(B^\top P+D^\top P\wt C+\wt S\big)\subseteq\sR\big(R+D^\top PD\big),\\
\ns\ds R+D^\top PD\ges0 \ea\right.$$
admits a (unique) static stabilizing solution $P\in\dbS^n$, and the following BSDE
\bel{11.7-2}\ba{ll}
\ds d\eta(t)=-\Big\{\big[\wt A^\top\1n-\big(PB+\wt C^\top PD+\wt S^\top\big)\big(R+D^\top PD\big)^\dag B^\top\big]\eta\\
\ns\ds\hphantom{d\eta(t)=-\Big\{}
+\big[\wt C^\top\1n-\big(PB+\wt C^\top PD+\wt S^\top\big)\big(R+D^\top PD\big)^\dag D^\top\big]\z\\
\ns\ds\hphantom{d\eta(t)=-\Big\{}
+\big[\wt C^\top\1n-\big(PB+\wt C^\top PD+\wt S^\top\big)\big(R+D^\top PD\big)^\dag D^\top\big]P\si\\
\ns\ds\hphantom{d\eta(t)=-\Big\{}
-\big(PB+\wt C^\top PD+\wt S^\top\big)\big(R+D^\top PD\big)^\dag\rho+Pb+\wt q\Big\}dt+\z dW(t), \qq t\ges0
\ea\ee
admits an $L^2$-stable adapted solution $(\eta(\cd),\z(\cd))$ such that
\bel{11.7-3}B^\top\eta(t)+D^\top\z(t)+D^\top P\si(t)+\rho(t)\in\sR\big(R+D^\top PD\big),\q\ae~t\in[0,\i),~~\as\ee
Using \rf{11.7-3} and the equality
$$\big(PB+\wt C^\top PD+\wt S^\top\big)\big(R+D^\top PD\big)^\dag
=\cL(P)\cN(P)^\dag+\Si^\top\cN(P)\cN(P)^\dag,$$
it is straightforward to show that \rf{11.7-2} is equivalent to
\begin{eqnarray*}
&& d\eta(t)=-\Big\{\big[A^\top\1n-\cL(P)\cN(P)^\dag B^\top\big]\eta+\big[C^\top\1n-\cL(P)\cN(P)^\dag D^\top\big]\z\\
&&\hphantom{d\eta(t)=-\Big\{}
+\big[C^\top\1n-\cL(P)\cN(P)^\dag D^\top\big]P\si-\cL(P)\cN(P)^\dag\rho+Pb+q\Big\}dt+\z dW(t), \qq t\ges0.
\end{eqnarray*}
Thus, we need only show that $P$ is a static stabilizing solution of ARE \rf{ARE}.
To this end, choose $\L\in\dbR^{m\times n}$ such that
$$\Si^*\deq-\big(R+D^\top PD\big)^\dag\big(B^\top P+D^\top P\wt C+\wt S\big)
+\big[I-\big(R+D^\top PD\big)^\dag\big(R+D^\top PD\big)\big]\L$$
is a stabilizer of $[\wt A,\wt C;B,D]$. We have
\bel{16Spet6-16:45}\ba{lll}
\ds\big(R+D^\top PD\big)(\Si^*+\Si)
\4n&=\4n&\ds -\,\big(B^\top P+D^\top P\wt C+\wt S\big)+\big(R+D^\top PD\big)\Si\\
\ns\4n&=\4n&\ds -\,\big(B^\top P+D^\top PC+S\big).
\ea\ee
It follows that $\sR\big(B^\top P+D^\top PC+S\big)\subseteq\sR\big(R+D^\top PD\big)$. Moreover,
\begin{eqnarray*}
0\4n&=\4n& P\wt A+\wt A^\top P+\wt C^\top P\wt C+\wt Q
-\big(PB+\wt C^\top PD+\wt S^\top\big)\big(R+D^\top PD\big)^\dag\big(B^\top P+D^\top P\wt C+\wt S\big)\\
\4n&=\4n& P(A+B\Si)+(A+B\Si)^\top P+(C+D\Si)^\top P(C+D\Si)+Q+S^\top\Si+\Si^\top S+\Si^\top R\Si\\
\4n&~\4n& -\,\big[PB+(C+D\Si)^\top PD+(S+R\Si)^\top\big]\big(R+D^\top PD\big)^\dag
\big[B^\top P+D^\top P(C+D\Si)+S+R\Si\big]\\
\4n&=\4n& PA+A^\top P+C^\top PC+Q
+\big(PB+C^\top PD+S^\top\big)\Si+\Si^\top\big(B^\top P+D^\top PC+S\big)\\
\4n&~\4n& -\,\big(PB+C^\top PD+S^\top\big)\big(R+D^\top PD\big)^\dag\big(B^\top P+D^\top PC+S\big)\\
\4n&~\4n& -\,\big(PB+C^\top PD+S^\top\big)\big(R+D^\top PD\big)^\dag\big(R+D^\top PD\big)\Si\\
\4n&~\4n& -\,\Si^\top\big(R+D^\top PD\big)\big(R+D^\top PD\big)^\dag\big(B^\top P+D^\top PC+S\big)\\
\4n&=\4n& PA+A^\top P+C^\top PC+Q
-\big(PB+C^\top PD+S^\top\big)\big(R+D^\top PD\big)^\dag\big(B^\top P+D^\top PC+S\big)\\
\4n&~\4n& +\,\big(PB+C^\top PD+S^\top\big)\big[I-\big(R+D^\top PD\big)^\dag\big(R+D^\top PD\big)\big]\Si\\
\4n&~\4n& +\,\Si^\top\big[I-\big(R+D^\top PD\big)\big(R+D^\top PD\big)^\dag\big]\big(B^\top P+D^\top PC+S\big)\\
\4n&=\4n& PA+A^\top P+C^\top PC+Q
-\big(PB+C^\top PD+S^\top\big)\big(R+D^\top PD\big)^\dag\big(B^\top P+D^\top PC+S\big)\\
\4n&~\4n& -\,(\Si^*+\Si)^\top\big(R+D^\top PD\big)\big[I-\big(R+D^\top PD\big)^\dag\big(R+D^\top PD\big)\big]\Si\\
\4n&~\4n& -\,\Si^\top\big[I-\big(R+D^\top PD\big)\big(R+D^\top PD\big)^\dag\big]\big(R+D^\top PD\big)(\Si^*+\Si)\\
\4n&=\4n& PA+A^\top P+C^\top PC+Q
-\big(PB+C^\top PD+S^\top\big)\big(R+D^\top PD\big)^\dag\big(B^\top P+D^\top PC+S\big).
\end{eqnarray*}
Therefore, $P$ solves ARE \rf{ARE}. Noting that $\Si^*+\Si$ is a stabilizer of $[A,C;B,D]$,
it is further clear from \rf{16Spet6-16:45} and Lemma \ref{Penrose} that $P$ is static stabilizing.
\end{proof}

\section{The One-Dimensional Case}

In this section, we look at the case where both the state and the control variables are one-dimensional.
For such a case, we can solve Problem (LQ)$^0$ completely. To avoid trivial exceptions we assume that
\bel{1d-1}\left\{\2n\ba{ll}
\ds B\neq0\hb{~~\,or\,~~}  D\neq0,\\
\ns\ds \sS[A,C;B,D]\neq\emptyset.\ea\right.\ee
By Lemma \ref{AC-stable}, the second condition in \rf{1d-1} is equivalent to the solvability of
$2(A+B\Th)+(C+D\Th)^2<0$ with the unknown $\Th$.
It is then easy to verify that \rf{1d-1} holds if and only if
\bel{1d-2}(2A+C^2)D^2<(B+CD)^2.\ee

\ss

Let us first look at the case $D=0$. By scaling, we may assume without loss of generality that $B=1$.
Then ARE \rf{ARE} becomes
\bel{ARE-D=0}\left\{\2n\ba{ll}
\ds P(2A+C^2)+Q-R^\dag(P+S)^2=0,\\
\ns\ds P+S=0\q\hb{if}~R=0,\\
\ns\ds R\ges0.\ea\right.\ee
Also, we note that, by Lemma \ref{AC-stable}, $\Th$ is a stabilizer of $[A,C;1,0]$ if and only if
$\Th<-(2A+C^2)/2$.

\bt{bt-D=0}\sl Suppose that $D=0$ and $B=1$. We have the following:\\[-1.6em]
\begin{enumerate}[\indent\rm(i)]
\item If $R<0$, then Problem {\rm(LQ)}$^0$ is not solvable.\\[-1.6em]
\item If $R=0$, then Problem {\rm(LQ)}$^0$ is solvable if and only if $Q=S(2A+C^2)$. In this case,
$$(\Th, v(\cd))\q\hb{with}\q \Th<-(2A+C^2)/2,\q v(\cd)\in L^2_\dbF(\dbR)$$
are all the closed-loop optimal strategies of Problem {\rm(LQ)}$^0$.\\[-1.6em]
\item If $R>0$, then Problem {\rm(LQ)}$^0$ is solvable if and only if
$\Si\equiv R(2A+C^2)^2-4S(2A+C^2)+4Q>0$. In this case,
$$\lt(-\,{2A+C^2+\sqrt{\Si/R}\over2},~0\rt)$$
is the unique closed-loop optimal strategy of Problem {\rm(LQ)}$^0$.
\end{enumerate}\et

\begin{proof}
(i) It is obvious because ARE \rf{ARE-D=0} is not solvable in this case.

\ms

(ii) When $R=0$, ARE \rf{ARE-D=0} further reduces to
$$\left\{\2n\ba{ll}
\ds P(2A+C^2)+Q=0,\\
\ns\ds P+S=0,\ea\right.$$
which is solvable if and only if $Q=S(2A+C^2)$. In this case, $\cN(P)=R=0$,
and the second assertion follows immediately from Corollary \ref{main-bc}.

\ms

(iii) When $R>0$, ARE \rf{ARE-D=0} can be written as
\bel{Ric-R>0}P^2+\big[2S-(2A+C^2)R\big]P+S^2-QR=0,\ee
which is solvable if and only if
$$0\les\D\deq\big [2S-(2A+C^2)R\big]^2-4(S^2-QR)=R\big[R(2A+C^2)^2-4S(2A+C^2)+4Q\big].$$
In the case of $\D\ges0$, \rf{Ric-R>0} has two solutions:
$$P_1={(2A+C^2)R-2S-\sqrt{\D}\over2},\qq P_2={(2A+C^2)R-2S+\sqrt{\D}\over2},$$
and $P_i$ is static stabilizing if and only if
$$-{2A+C^2\over2}>-{P_i+S\over R}=-{2A+C^2\over2}-{(-1)^i\sqrt{\D}\over2R}.$$
Clearly, $P_1$ cannot be static stabilizing, and $P_2$ is static stabilizing if and only if $\D>0$,
or equivalently, $R(2A+C^2)^2-4S(2A+C^2)+4Q>0$. The second assertion follows easily.
\end{proof}

We now look at the case $D\neq0$. \rm As before, we may assume, without loss of generality
(by scaling, if necessary), that $D=1$. Denote
\bel{abg}\left\{\2n\ba{ll}
\ds \a=(B+C)^2-(2A+C^2),\\
\ns\ds \b=Q-(2A+C^2)R+2(B+C)[(B+C)R-S],\\
\ns\ds \g=[(B+C)R-S]^2.\ea\right.\ee
Then \rf{1d-2} is equivalent to $\a>0$, and $\Th$ is a stabilizer of $[A,C;B,1]$ if and only if
\bel{Th-D=1}|\Th+B+C|<\sqrt{\a}.\ee

\bt{bt-D=1}\sl Suppose that $D=1$ and $\a>0$. Then Problem {\rm(LQ)}$^0$ is solvable
if and only if one of the following conditions holds:

\ms

{\rm(i)} $Q=(2A+C^2)R$ and $S=(B+C)R$. In this case,
\bel{R+P=0}(\Th, v(\cd))\q\hb{with}\q |\Th+B+C|<\sqrt{\a},\q v(\cd)\in L^2_\dbF(\dbR)\ee
are all the closed-loop optimal strategies of Problem {\rm(LQ)}$^0$.

\ms

{\rm(ii)} $2A+C^2\neq0$, $(2A+C^2)S\ges(B+C)Q$, and
$$R>{2(B+C-\sqrt{\a})S-Q\over(B+C-\sqrt{\a})^2}.$$

\ss

{\rm(iii)} $2A+C^2\neq0$, $(2A+C^2)S<(B+C)Q$, and
$$R>{2(B+C+\sqrt{\a})S-Q\over(B+C+\sqrt{\a})^2}.$$

\ss

{\rm(iv)} $2A+C^2=0$, $Q>0$, and
$$R>{4(B+C)S-Q\over4(B+C)^2}.$$

\ss

In the cases {\rm(ii)}, {\rm(iii)}, and {\rm(iv)},
$$\lt({2\a[(B+C)R-S]\over \b+\sqrt{\b^2-4\a\g}}-(B+C),~0\rt)$$
is the unique closed-loop optimal strategy of Problem {\rm(LQ)}$^0$.
\et

\begin{proof}
We rewrite ARE \rf{ARE} as follows:
\bel{D=1}\left\{\2n\ba{ll}
\ds P(2A+C^2)+Q-(R+P)^\dag[P(B+C)+S]^2=0,\\
\ns\ds P(B+C)+S=0\q\hb{if}~R+P=0,\\
\ns\ds R+P\ges0.\ea\right.\ee
By Corollary \ref{main-bc}, Problem {\rm(LQ)}$^0$ is solvable if and only if
\rf{D=1} admits a static stabilizing solution.
So we need only discuss the solvability of \rf{D=1}.

\ms

Clearly, $P=-R$ is a solution of \rf{D=1} if and only if
\bel{QSR}Q=(2A+C^2)R,\qq S=(B+C)R.\ee
In this case, $P=-R$ is also static stabilizing, and $\cN(P)=R+P=0$.
By Corollary \ref{main-bc} and \rf{Th-D=1},
we see all closed-loop optimal strategies of Problem {\rm(LQ)}$^0$ are given by \rf{R+P=0}.

\ms

If \rf{QSR} does not hold, by the change of variable $y=R+P$,
equation \rf{D=1} further reduces to the following:
\bel{D=1-y}\left\{\2n\ba{ll}
\ds \a y^2-\b y+\g=0,\\
\ns\ds y>0,\ea\right.\ee
which is solvable if and only if
$$\D=\b^2-4\a\g\ges0,\qq \b+\sqrt{\D}>0,$$
or equivalently (noting that $\a>0$ and $\g\ges0$),
\bel{}\D=\b^2-4\a\g\ges0,\qq \b>0.\ee
In this case, if $\g>0$, then \rf{D=1-y} has two solutions:
\bel{g>0_y12}y_1=R+P_1={\b-\sqrt{\D}\over2\a},\qq y_2=R+P_2={\b+\sqrt{\D}\over2\a}.\ee
For $i=1,2$, let
$$\Th_i=-{P_i(B+C)+S\over R+P_i}={(B+C)R-S\over y_i}-(B+C).$$
Note that $\Th_i$ is a stabilizer of $[A,C;B,1]$ if and only if $|\Th_i+B+C|<\sqrt{\a}$,
or equivalently,
\bel{16Step7-18:00}\g=[(B+C)R-S]^2<\a y_i^2=\b y_i-\g.\ee
Upon substitution of \rf{g>0_y12} into \rf{16Step7-18:00}, the latter is in turn equivalent to
\bel{11.23-1}\D+(-1)^i\b\sqrt{\D}>0.\ee
Obviously, \rf{11.23-1} cannot hold for $i=1$, and it holds for $i=2$ if and only if $\D>0$.
Likewise, if $\g=0$, then $P_2$ is the unique solution of \rf{D=1-y},
and $\Th_2$ is a a stabilizer of $[A,C;B,1]$ if and only if $\D>0$.
Therefore, ARE \rf{D=1} admits a static stabilizing solution $P\neq R$ if and only if
\bel{ii-iii}\b>0,\qq \b^2-4\a\g>0.\ee
Noting that $\b=[(B+C)^2+\a]R+Q-2(B+C)S$, we have by a straightforward computation:
\begin{eqnarray*}
\b^2-4\a\g \4n&=\4n& [(B+C)^2-\a]^2R^2-\(4[(B+C)^2-\a](B+C)S-2[(B+C)^2+\a]Q\)R\\
\4n&~\4n& +\,Q^2-4(B+C)SQ+4[(B+C)^2-\a]S^2\\
\4n&\equiv\4n& aR^2-bR+c.
\end{eqnarray*}
Also, we have
$$b^2-4ac=16\a\([(B+C)^2-\a]S-(B+C)Q\)^2\ges0.$$
If $a=[(B+C)^2-\a]^2=(2A+C^2)^2\neq0$, then $\b^2-4\a\g>0$ if and only if
$$ R>{b+\sqrt{b^2-4ac}\over2a},\q\hb{or}\q R<{b-\sqrt{b^2-4ac}\over2a}.$$
Because $R={2(B+C)S-Q\over(B+C)^2+\a}$ implies $\b^2-4\a\g=-4\a\g\les0$, we have
$${b-\sqrt{b^2-4ac}\over2a}\les{2(B+C)S-Q\over(B+C)^2+\a}\les{b+\sqrt{b^2-4ac}\over2a}.$$
Therefore, \rf{ii-iii} holds if and only if
\bel{}R>{b+\sqrt{b^2-4ac}\over2a}=\left\{\2n\ba{ll}
\ds {2(B+C-\sqrt{\a})S-Q\over(B+C-\sqrt{\a})^2},\qq \hb{if}\q(2A+C^2)S\ges(B+C)Q,\\
\ns\ds {2(B+C+\sqrt{\a})S-Q\over(B+C+\sqrt{\a})^2},\qq \hb{if}\q(2A+C^2)S<(B+C)Q.
\ea\right.\ee
If $a=[(B+C)^2-\a]^2=(2A+C^2)^2=0$, then
\begin{eqnarray*}
\b \4n&=&\4n 2(B+C)^2R-2(B+C)S+Q,\\
\b^2-4\a\g \4n&=&\4n Q\big[4(B+C)^2R-4(B+C)S+Q\big],
\end{eqnarray*}
and it is not hard to see that \rf{ii-iii} holds if and only if
\bel{}Q>0,\qq R>{4(B+C)S-Q\over4(B+C)^2}.\ee

\ss

Finally, in the cases (ii), (iii), and (iv), we see from the preceding argument that
ARE \rf{D=1} has a unique static stabilizing solution
$$P={\b+\sqrt{\D}\over2\a}-R.$$
Note that $\cN(P)=R+P>0$ and
$$-\cN(P)^{-1}\cL(P)=-{P(B+C)+S\over R+P}
={2\a[(B+C)R-S]\over \b+\sqrt{\b^2-4\a\g}}-(B+C).$$
The last assertion follows immediately from Corollary \ref{main-bc}.
\end{proof}


\end{document}